\documentclass[12pt, reqno]{amsart} 

\usepackage{xcolor} 
\usepackage[colorlinks=true,linkcolor=blue,citecolor=red,urlcolor=red,backref=page]{hyperref} 
\usepackage{esint} 
\usepackage{amssymb} 
\usepackage{enumitem} 
\usepackage{mathrsfs}

\setenumerate{label=(\alph*)} 
\setenumerate[2]{label=(\Alph*)} 

\makeatletter 
\def\l@subsection{\@tocline{2}{0pt}{2.8pc}{5pc}{}}

\newtheorem{theorem}{Theorem}[section]
\newtheorem{proposition}[theorem]{Proposition}
\newtheorem{lemma}[theorem]{Lemma}
\newtheorem{corollary}[theorem]{Corollary}

\theoremstyle{definition}
\newtheorem{definition}[theorem]{Definition}

\theoremstyle{remark}
\newtheorem{remark}[theorem]{Remark}

\numberwithin{equation}{section} 

\def\div{\operatorname{div}}
\def\dist{\operatorname{dist}}

\def\Ker{\operatorname{Ker}}
\def\Im{\operatorname{Im}}
\def\Coker{\operatorname{Coker}}
\def\ind{\operatorname{ind}}
\def\BMO{\operatorname{BMO}}
\def\VMO{\operatorname{VMO}}

\def\ess{\operatorname{ess}}
\def\inf{\operatornamewithlimits{inf\vphantom{p}}}

\newcommand{\norm}[1]{\left\lVert#1\right\rVert}

\DeclareMathOperator*{\esssup}{ess\,sup}
\DeclareMathOperator*{\essinf}{ess\,inf}

\calclayout

\begin{document}

\title[The BA property and interpolation of the essential norm]{The bounded approximation property in Banach function lattices, interpolation of the essential norm, and applications}

\author{Juan Jos\'{e} Mar\'{i}n}
\address{Juan Jos\'{e} Mar\'{i}n
\\
Departamento de Matem{\'a}tica Aplicada y Estad{\'i}stica
\\
Universidad Polit{\'e}cnica de Cartagena
\\
Cartagena, Spain.} \email{juanjose.marin@upct.es}




\begin{abstract}
We prove the bounded approximation property for Banach function lattices with order continuous norm on which the Hardy–Littlewood maximal operator is bounded.
This allows us to develop a general framework for the quantitative interpolation of compact operators, understood in terms of the essential norm, encompassing a large class of function spaces including weighted Lebesgue spaces, rearrangement invariant Banach function spaces, and variable exponent Lebesgue spaces.
Combined with extrapolation techniques, we show how to transfer compactness  (or, more generally, control of the essential norm) from unweighted Lebesgue spaces to the three aforementioned classes of function spaces.
Finally, we present several applications in analysis and partial differential equations.
\end{abstract}

\maketitle

\tableofcontents

\section{Introduction}

Interpolation theory provides a powerful framework for transferring properties of operators between Banach spaces, yet the behavior of compactness under interpolation remains delicate. Specifically, if a linear operator acts compactly on one endpoint space of a compatible couple and boundedly on the other, it is natural to ask whether compactness persists in the intermediate spaces generated by interpolation methods. 
This problem can be traced back to the work of Krasnosel'ski{\u\i} \cite{Kra60} in 1960,  where the case of Lebesgue spaces is addressed,
and was subsequently developed for more general Banach spaces in the seminal works of
Lions and Peetre \cite{LioPee64} and Calderón \cite{Cal64} in 1964.
While this question was answered affirmatively for the real interpolation method in 1992 by Cwikel \cite{Cwi92} and Cobos, Kühn, and Schonbek \cite{CobKuhSch92}, the problem for the complex interpolation method in full generality is still open.
Since then, numerous authors have contributed to this topic. In particular, we highlight two works closely related to the present paper, as they address the interpolation of compact operators in Banach lattices via general interpolation methods: the work of Pustylnik \cite{Pus08} establishes results in the setting of pseudo-lattices of measurable functions, while Brudnyi \cite{Bru15} develops the theory for abstract Banach lattices.

Quantitative interpolation of compact operators refines the previous question by seeking explicit control on the degree of compactness across intermediate spaces. Rather than merely asserting that an operator remains compact, this approach measures compactness via tools such as the measure of non-compactness or the essential norm (cf., e.g., \cite{AksTyl05,CobFerMar99,Fer06}).
We are interested in quantifying how the essential norm behaves under interpolation.

The essential norm of a linear bounded operator $T:X\to Y$ between two Banach spaces $X$ and $Y$ is defined as
\begin{equation}\label{eq:ess-def}
\|T\|^{\ess}_{X\to Y} := \dist(T,\mathcal{K}(X\to Y)) = \inf_{K\in\mathcal{K}(X\to Y)} \|T-K\|_{X\to Y},
\end{equation}
where $\mathcal{K}(X\to Y)$ denotes the space of all compact operators $K:X\to Y$.
Assume that $(X,\mathbb{Y})$ is an interpolation pair between two compatible couples of Banach spaces $(X_0,X_1)$ and $(\mathbb{Y}_0, \mathbb{Y}_1)$, in the sense that there is a continuous function $\Phi$ which is non-decreasing in each variable such that
\begin{equation}
\|T\|_{X\to \mathbb{Y}} \leq \Phi \left(\|T\|_{X_0\to \mathbb{Y}_0}, \, \|T\|_{X_1 \to \mathbb{Y}_1}\right),
\end{equation}
for any linear operator $T:X_0+X_1\to \mathbb{Y}_0+\mathbb{Y}_1$ with bounded restrictions $T:X_0\to \mathbb{Y}_0$ and $T:X_1\to \mathbb{Y}_1$. In this context, we seek to identify the set of assumptions under which one has
\begin{equation}\label{eq:intro.ess}
\|T\|^{\ess}_{X\to \mathbb{Y}} \leq \Phi \left(C_0 \|T\|^{\ess}_{X_0\to \mathbb{Y}_0}, \, C_1 \|T\|^{\ess}_{X_1 \to \mathbb{Y}_1}\right).
\end{equation}
Theorem \ref{thm:general} addresses this question assuming that $\mathbb{Y}_0$, $\mathbb{Y}_1$, and $\mathbb{Y}$ are Banach function lattices over a  space of homogeneous type, the Hardy–Littlewood maximal operator is bounded on each of them, and both $\mathbb{Y}_0$ and $\mathbb{Y}_1$ possess order continuous norms (see Sections \ref{sec:prelim} and \ref{sec:BL} for the definitions).
This may be regarded as a quantitative extension of the work of Pustylnik in \cite{Pus08}.
It also complements the results of Aksoy and Tylli \cite{AksTyl05}, established for the real interpolation method under different assumptions.

The proof of the aforementioned result relies on the bounded approximation property, which asserts that the identity operator can be uniformly approximated on compact subsets by finite rank bounded linear operators whose norms are uniformly bounded (cf. Definition \ref{def:BAP}).
The bounded approximation property bridges finite-dimensional and infinite-dimensional analysis and constitutes a fundamental concept for the theory of Banach spaces, underpinning many results in functional analysis.
Specifically, we show (cf. Theorem \ref{thm:BAP-Xw}) that 
\begin{equation}
\parbox{11.9cm}{\begin{center}
any Banach function lattice over a space of homogeneous type, with order continuous norm and bounded Hardy–Littlewood maximal operator, has the bounded approximation property. 
\end{center}}
\end{equation}
Beyond its role as a key ingredient in our analysis of compact interpolation, this result is of intrinsic interest. In particular, it provides new insight into the bounded approximation property in (weighted) rearrangement invariant Banach function spaces and in (weighted) variable exponent Lebesgue spaces.
An additional benefit of the techniques developed in this work, at the expense of specializing to Banach function lattices on spaces of homogeneous type, is the explicit nature of the construction of the approximating operators.

We emphasize that, in contrast to the results in \cite{AksTyl05} for the real interpolation method, our interpolation theorem does not require the intermediate space $\mathbb{Y}$ to have the bounded approximation property. This overcomes a limitation, since this property is not preserved by interpolation and hence the required condition on $\mathbb{Y}$ cannot be inferred from the same property on the endpoint spaces.

In Section \ref{sec:interpol}, we specialize our interpolation result to three prominent classes of function spaces: weighted Lebesgue spaces, weighted rearrangement invariant Banach function spaces (such as Orlicz spaces), and variable exponent Lebesgue spaces. We show that all these cases fall within our general interpolation framework and combine this with extrapolation techniques to transfer control of the essential norm (and, in particular, compactness) from unweighted Lebesgue spaces to each of the function spaces mentioned above.

Finally, Section \ref{sec:applications} focuses on applications to analysis and partial differential equations.  While some of the results in this section are essentially known (cf. \cite{MitMitMit-book3, MitMitMit-book5}), our approach introduces original techniques and provides more direct proofs.
The first application addresses bounds for the essential norm of commutators $[M_b,T]$, defined as $[M_b,T]f = bT(f) - T(bf)$, with $T$ a Calderón–Zygmund operator and $b \in \mathrm{BMO}(\Sigma,\mu)$. This generalizes the classical result of Uchiyama (cf. \cite{Uch78}), which asserts that $[M_b,T]$ is compact in Lebesgue spaces if and only if $b \in \mathrm{VMO}(\Sigma,\mu)$, extending it to a broader class of function spaces in a quantitative manner.

We next turn to applications in Fredholm theory. It is well known that the operator $\lambda I+T$ is invertible via a Neumann series whenever $T$ has sufficiently small norm relative to $\lambda \in \mathbb{R}\setminus\{0\}$. Similarly, $\lambda I+T$ is Fredholm of index zero whenever $T$ has sufficiently small essential norm. In this context, our transfer techniques allow us to establish Fredholmness in a wide range of function spaces under the assumption that the essential norm is small only in classical Lebesgue spaces.

Our third application concerns estimates for the essential norm of a class of singular integral operators defined on the boundary of a rough domain, whose kernels involve an inner product with the outward unit normal vector. Notable examples include the harmonic double layer potential operators.
When combined with the Fredholm theory developed earlier, these results yield invertibility of layer potential operators. This provides the basis for our final application, which addresses the transfer of well-posedness for the Dirichlet problem in non-smooth domains from Lebesgue spaces to a variety of other function spaces.

\section{Preliminary matters} \label{sec:prelim}

We start by defining spaces of homogeneous type, introduced by Coifman and Weiss in \cite{CoiWei71}. Recall that a {\tt quasi-distance} $d$ on a set $\Sigma$ is a function satisfying the same conditions as a distance, except that the triangle inequality is weakened to 
\begin{equation} \label{eq:triangularineq}
d(x,y) \leq A_0 (d(x,z)+d(z,y)) \quad \text{ for every } x,y,z\in\Sigma,
\end{equation}
where $A_0\in [1,\infty)$.

\begin{definition}
A {\tt space of homogeneous type} is a triple $(\Sigma,d,\mu)$, where $\Sigma$ is a non-empty set 
endowed with a quasi-distance $d$ for which balls
\begin{equation}
B(x,r) := \{y \in \Sigma : d(x,y)<r \}
\text{ are open for every } x \in \Sigma \text{ and } r>0,
\end{equation}
and $\mu$ is a non-negative regular Borel measure satisfying 
the doubling condition, i.e., there is a constant $A_1\in [1,\infty)$ such that
\begin{equation}\label{eq:doubling}
0<\mu(B(x,2r)) \leq A_1 \mu(B(x,r))<\infty \quad \text{ for every } x \in \Sigma \text{ and } r>0.
\end{equation}
\end{definition}
Note that every space of homogeneous type is $\sigma$-finite, since
$\Sigma=\cup_{n=1}^{\infty}B(x_0,n)$ for any $x_0\in\Sigma$ and we are assuming that balls have finite measure.
As balls also have positive measure, we have $\mu(\Sigma)>0$ as well.
We will sometimes assume that $(\Sigma, \mu)$ is {\tt non-atomic}, that is,
for every $\mu$-measurable set $E$ with $\mu(E)>0$, there is a $\mu$-measurable subset $F\subset E$ such that $0<\mu(F)<\mu(E)$.
We say that $(\Sigma,\mu)$ is {\tt complete} if every subset of each set of measure zero is measurable.

For every function $f$ which is locally integrable in $\Sigma$, we define the {\tt Hardy-Littlewood} {\tt maximal} {\tt operator} as
\begin{equation}\label{eq:maximal}
\mathcal{M} f(x):=\sup_{B \ni x}\fint_{B}|f(y)|\,d\mu(y),\qquad 
\text{ for each }\,x\in\Sigma,
\end{equation}
where the supremum is taken over all balls $B$ in $\Sigma$ containing
$x$ and the barred integral represents the average.

We continue by considering the {\tt dyadic} {\tt decomposition} of M. Christ (cf. \cite{Chr90b}
and the refinement in \cite{HytKai12}) in the context of spaces of homogeneous type.

\begin{theorem}\label{thm:dyadic-dec}
Let $(\Sigma,d,\mu)$ be a space of homogeneous type. Then, there exists a collection of sets 
\begin{equation} \label{eq:dyadic-dec}
\mathbb{D}(\Sigma) =
\{Q_{\alpha}^k \subseteq \Sigma : k \in \mathbb{Z}, \enskip \alpha\in I_k \},
\end{equation}
where $I_k$ is a countable index set,
that satisfies the following properties:
\begin{enumerate}
\item For every $k \in \mathbb{Z}$,
\begin{equation}
\Sigma = \bigcup_{\alpha\in I_k} Q_{\alpha}^k.
\end{equation}
\item If $Q_\alpha^{k},  Q_\beta^{\ell} \in \mathbb{D}(\Sigma)$, then either $Q_\alpha^{k} \cap Q_\beta^{\ell} = \varnothing$,
$Q_\alpha^{k} \subseteq Q_\beta^{\ell}$, or 
$Q_\beta^{\ell} \subseteq Q_\alpha^{k}$.
\item For any $Q^k_\alpha \in \mathbb{D}(\Sigma)$ there is at least one set
$Q^{k+1}_\beta \in \mathbb{D}(\Sigma)$ and exactly one set $Q^{k-1}_\gamma \in \mathbb{D}(\Sigma)$ such that
\begin{equation}
Q^{k+1}_\beta \subset Q^k_\alpha\subset Q^{k-1}_\gamma.
\end{equation}
\item There is $C_0\in(1,\infty)$ such that 
$\mu(Q^k_\alpha) \leq C_0 \mu(Q^{k+1}_\beta)$
whenever $Q^{k+1}_\beta \subseteq Q^k_\alpha$.
\item There are $C_1,C_2\in(0,\infty)$ and $\delta\in(0,1)$ such that every $Q^k_\alpha \in \mathbb{D}(\Sigma)$ contains a ball of radius $C_1 \delta^k$ and is contained in a ball of radius $C_2 \delta^k$. 
\end{enumerate}
\end{theorem}

Next, define the dyadic maximal operator
\begin{equation}\label{eq:maximal-D}
\mathcal{M}_{\mathbb{D}} f(x):=\sup_{Q \ni x, \, Q \in \mathbb{D}(\Sigma)}\fint_{Q}|f(y)|\,d\mu(y),\qquad 
\text{ for each }\,x\in\Sigma.
\end{equation}
As a consequence of item $(e)$ above and the doubling condition, one can prove (cf. \cite[Corollary 7.4]{HytKai12}) that
there is a constant $C_\Sigma\in[1,\infty)$ such that for every
$Q \in \mathbb{D}(\Sigma)$ there is a ball $B_Q$ with $Q \subseteq B_Q$ and $\mu(B_Q)\leq C_\Sigma\mu(Q)$. Therefore, for every locally integrable function $f$ in $\Sigma$, we have
\begin{equation}\label{eq:maximal-compar}
\mathcal{M}_{\mathbb{D}} f(x) \leq \ C_\Sigma \mathcal{M} f(x), \qquad \forall x \in \Sigma.
\end{equation}

To close this section, we present essential definitions and properties concerning weights. A more detailed exposition of weight theory on spaces of homogeneous type can be found in \cite{GenGogKokKrb97,StrTor89}.
A {\tt weight} $w$ on $(\Sigma,\mu)$ is a $\mu$-measurable function such that $0<w(x)<\infty$ 
for $\mu$-a.e. point $x\in\Sigma$. 
Given a space of homogeneous type $(\Sigma,d,\mu)$ and
$p\in(1,\infty)$, we say that $w$ lies in the
{\tt Muckenhoupt class} $A_p(\Sigma)$ if 
\begin{equation}\label{eq:def-Apweight}
[w]_{A_p} := \sup_{B \subseteq \Sigma} \left( \fint_{B} w(x) d\mu (x) \right) 
\left( \fint_{B} w(x)^{1-p'} d\mu (x) \right) ^{p-1} < \infty,
\end{equation}
where the supremum is taken over all balls
$B \subseteq \Sigma$ and $p'$ is the conjugate exponent of $p$, that is, $1/p+1/p'=1$. Similarly, $w\in A_1(\Sigma)$ if 
\begin{equation}\label{eq:def-A1weight}
[w]_{A_1} := \sup_{B \subseteq \Sigma} \left( \essinf_{x \in B} w(x) \right)^{-1} 
\left( \fint_{B} w(x) d\mu (x) \right) < \infty.
\end{equation}
In particular, 
\begin{equation}
\mathcal{M}w(x)\leq[w]_{A_1}w(x) \qquad \text{for }\mu\text{-a.e. }x\in\Sigma.
\end{equation}
We also define 
\begin{equation}
A_{\infty}(\Sigma) :=
\bigcup_{p \geq 1} A_p(\Sigma).
\end{equation}
From \eqref{eq:def-Apweight}, one readily observes that for every $p\in(1,\infty)$,
\begin{equation}\label{eq:Ap-p1}
\parbox{7.5cm}{\begin{center}
$w \in A_p(\Sigma)$ if and only if $w^{1-p'}\in A_{p'}(\Sigma)$ and in that case, $[w^{1-p'}]_{A_{p'}} = [w]_{A_p}^{1/(p-1)}$
\end{center}}
\end{equation}
and, by Hölder’s inequality, we obtain that
\begin{equation}\label{eq:Ap-p2}
\parbox{14cm}{\begin{center}
if $1\leq p\leq q\leq \infty$ then $A_p(\Sigma) \subseteq A_q(\Sigma)$.
\end{center}}
\end{equation}
Moreover, one has the following openness property:
\begin{equation}\label{eq:openness}
\parbox{13cm}{\begin{center}
if $w \in A_p(\Sigma)$ with $p\in(1,\infty)$ then there exists $\tau\in(1,\infty)$ such that $w^{\tau} \in A_{p}(\Sigma)$ and there exists $\varepsilon\in(0,p-1)$ such that $w \in A_{p-\varepsilon}(\Sigma)$.
\end{center}}
\end{equation}
Throughout, $w$ will represent the weighted measure $w \mu$, that is, for every $\mu$-measurable set $E\subseteq\Sigma$ we define 
\begin{equation}
w(E):=\int_E w \, d\mu.
\end{equation}
In this context, it is well-known (see, e.g., \cite[Proposition 7.13]{HytKai12}) that:
\begin{equation}\label{eq:Mbounded-Lpw}
\parbox{13cm}{\begin{center}
for any $p\in(1,\infty)$, ${\mathcal{M}}$ is bounded on $L^p(\Sigma,w)$ if and only if $w \in A_p(\Sigma)$ and in that case, $\|\mathcal{M}\|_{L^p(\Sigma,w) \to L^p(\Sigma,w)} \leq C_{\Sigma,\mu,p} [w]_{A_p}^{1/(p-1)}$.
\end{center}}
\end{equation}
%

\section{Banach function lattices and spaces}\label{sec:BL}

We now introduce Banach function lattices within the framework of a general $\sigma$-finite measure space $(\Sigma,\mu)$.

\begin{definition}\label{def:BL}
Let $(\Sigma,\mu)$ be a $\sigma$-finite measure space and let $\mathbb{M}_{\mu}$ be the set of all real-valued
$\mu$-measurable functions on $\Sigma$. A Banach space $(\mathbb{X},||\!\cdot\!||_{\mathbb{X}})$ is called {\tt Banach} {\tt function} {\tt lattice} over $(\Sigma,\mu)$ if $\mathbb{X}\subseteq\mathbb{M}_{\mu}$ and
\begin{equation}\label{eq:BL}
||f||_{\mathbb{X}} \leq ||g||_{\mathbb{X}} \quad
\text{whenever} \,\,|f|\leq|g| \,\, \mu\text{-a.e. on }\Sigma.
\end{equation}
\end{definition}

Throughout this work, we agree to use blackboard bold font whenever a Banach space is a function lattice.

\begin{definition}\label{def:OI}
Let $(\Sigma,\mu)$ be a $\sigma$-finite measure space and
$(\mathbb{X},\|\!\cdot\!\|_{\mathbb{X}})$ a Banach function lattice over $(\Sigma,\mu)$.
We say that $\mathbb{X}$ has {\tt order continuous norm} if for every $f\in\mathbb{X}$ and every decreasing sequence $\{ f_n \}_n \subseteq \mathbb{M}_\mu$ satisfying 
\begin{equation}\label{eq:ocn-condition}
0 \leq f_n \leq |f| \text{ for every } n\geq 1 \quad\text{ and }\quad f_n\to 0 \,\, \mu\text{-a.e. on }\Sigma \text{ as } n\to\infty,
\end{equation}
one has that 
\begin{equation}
\|f_n \|_{\mathbb{X}}\to 0 \,\, \text{ as } \,\, n\to\infty.
\end{equation}
\end{definition}

In view of later applications, we introduce the following weighted version of a Banach function lattice.

\begin{definition}
Let $(\Sigma,\mu)$ be a $\sigma$-finite measure space, $\mathbb{X}$ a Banach function lattice over $(\Sigma,\mu)$, and $w$ a weight on $(\Sigma,\mu)$. We define
\begin{equation}
\mathbb{X}_w = \{ f \in \mathbb{M} _{\mu} : \norm{f}_{\mathbb{X}_w} < \infty \} ,
\end{equation}
where $\|f\|_{\mathbb{X}_w} := \|fw\|_{\mathbb{X}}.$
\end{definition}

The following result guarantees that the spaces $\mathbb{X}_w$ defined above are Banach function lattices, and that the order continuity of their norm is equivalent to that of $\mathbb{X}$.
Even though the proof is elementary, we include it here for completeness.

\begin{proposition}\label{prop:Xw-properties}
Let $(\Sigma,\mu)$ be a $\sigma$-finite measure space, $\mathbb{X}$  a Banach function lattice over $(\Sigma,\mu)$ and $w$ a weight on $(\Sigma,\mu)$. Then,
\begin{enumerate}
\item $\mathbb{X}_w$ is a Banach function lattice over $(\Sigma,\mu)$.
\item $\mathbb{X}_w$ has order continuous norm if and only if $\mathbb{X}$ does.
\end{enumerate}
\end{proposition}
\begin{proof}

For $(a)$, it is easy to check that $\mathbb{X}_w$ is a normed space using that $w>0$ $\mu$-a.e. and the fact that $\mathbb{X}$ is a normed space.
If $\{f_n\}_n$ is a Cauchy sequence in $\mathbb{X}_w$ then 
$\{f_n w\}_n$ is a Cauchy sequence in $\mathbb{X}$. Since $\mathbb{X}$ is complete,
there is $g\in\mathbb{X}$ such that $f_n w \to g$ in $\mathbb{X}$. As $w>0$ $\mu$-a.e., we can define $f:=g/w$ with $fw=g\in\mathbb{X}$, that is, $f\in\mathbb{X}_w$. Hence,
\begin{equation}
\|f_n - f\|_{\mathbb{X}_w} = \|f_n w - g\|_{\mathbb{X}} \to 0 \quad \text{ as } n \to \infty.
\end{equation}
This proves that $\mathbb{X}_w$ is a Banach space. As $w>0$ $\mu$-a.e., the lattice property \eqref{eq:BL} for $\mathbb{X}_w$ follows at once from the lattice property for $\mathbb{X}$.

We now turn to the proof of $(b)$.
Assume that $\mathbb{X}$ has order continuous norm and let $f \in \mathbb{X}_w$ and $\{ f_n \}_n \subseteq \mathbb{M}_\mu$ a decreasing sequence satisfying \eqref{eq:ocn-condition}. As $0<w<\infty$ $\mu$-a.e., we have that $fw \in \mathbb{X}$ and $\{ f_n w \}_n \subseteq \mathbb{M}_\mu$ is a decreasing sequence satisfying
\begin{equation}
0 \leq f_n w \leq |fw| \text{ for every } n\geq 1 \quad\text{ and }\quad f_nw\to 0 \,\, \mu\text{-a.e. on }\Sigma \text{ as } n\to\infty.
\end{equation}
Using that $\mathbb{X}$ has order continuous norm, we finally obtain 
\begin{equation}
\|f_n \|_{\mathbb{X}_w}=\|f_n w \|_{\mathbb{X}}\to 0 \,\, \text{ as } \,\, n\to\infty.
\end{equation}
This proves that if $\mathbb{X}$ has order continuous norm, then $\mathbb{X}_w$ does as well. The converse implication also follows from this, since
$\mathbb{X}=(\mathbb{X}_w)_{1/w}$.
\end{proof}

One particularly relevant class of Banach function lattices is the class of Banach function spaces defined below.

\begin{definition}\label{def:BFS}
Let $(\Sigma,\mu)$ be a $\sigma$-finite measure space and let $\mathbb{M}_{\mu}$ be the set of all 
$\mu$-measurable functions. Given a mapping $\norm{\cdot}_{\mathbb{X}}:\mathbb{M}_{\mu} \to [0,\infty]$,
the set 
\begin{equation} 
\mathbb{X} = \{ f \in \mathbb{M} _{\mu} : \norm{f}_{\mathbb{X}} < \infty \} 
\end{equation} 
is called {\tt Banach} {\tt function} {\tt space} over $(\Sigma,\mu)$ if $\norm{\cdot}_{\mathbb{X}}$ satisfies the following 
properties for all $f,g \in \mathbb{M} _{\mu}$:
\begin{enumerate}
\item $\norm{f}_{\mathbb{X}} = \norm{ \left| f \right| }_{\mathbb{X}}$ and $\norm{f}_{\mathbb{X}}=0$ 
if and only if $f=0$ $\mu$-a.e. on $\Sigma$.
\vspace*{0.15cm}
\item $\norm{f+g}_{\mathbb{X}} \leq \norm{f}_{\mathbb{X}} + \norm{g}_{\mathbb{X}}$.
\vspace*{0.15cm}
\item $\norm{\lambda f}_{\mathbb{X}} = |\lambda| \norm{f}_{\mathbb{X}}$ for all $\lambda \in \mathbb{R}$.
\vspace*{0.15cm}
\item If $|f| \leq |g|$ $\mu$-a.e. on $\Sigma$, then $\norm{f}_{\mathbb{X}} \leq \norm{g}_{\mathbb{X}}$.
\vspace*{0.15cm}
\item If $\{ f_n \}_n \subseteq \mathbb{M}_\mu$ is a sequence such that $|f_n|$ increases to 
$|f|$ $\mu$-a.e. on $\Sigma$ as $n\to\infty$, then $\norm{f_n}_{\mathbb{X}}$ increases to $\norm{f}_{\mathbb{X}}$ as $n\to\infty$.
\vspace*{0.15cm}
\item If $E \subseteq \Sigma$ is a $\mu$-measurable set such that $\mu(E) < \infty$, then its
characteristic function $\mathbf{1}_E$ satisfies
$\norm{\mathbf{1}_E}_{\mathbb{X}} < \infty$ and there is $C_E\in(0,\infty)$ such that
$\int_E |f(x)| d\mu \leq C_E \norm{f}_{\mathbb{X}}$.
\end{enumerate}
\end{definition}

One can check that if the mapping $\norm{\cdot}_{\mathbb{X}}$
satisfies the conditions above, then $(\mathbb{X},\norm{\cdot}_{\mathbb{X}})$ is a complete normed space in accordance with the terminology.
As a consequence of this and properties $(d)$ and $(f)$ in the previous definition, every Banach function space is a Banach function lattice whose functions are locally integrable.

Given a Banach function space $\mathbb{X}$, we can define the associate space of 
$\mathbb{X}$, denoted by $\mathbb{X}'$, as the Banach function space defined through the norm
\begin{equation}
\norm{f}_{ \mathbb{X}^{'}} = \sup \left\{ \int_{\Sigma} |f(x)g(x)| d\mu(x) : 
g \in  \mathbb{X}, \enskip \norm{g}_{\mathbb{X}} \leq 1 \right\}.
\end{equation}
As a consequence, one obtains the following generalized H\"{o}lder inequality:
\begin{equation} \label{eq:XHolder}
\int_{\Sigma} |f(x)g(x)| d\mu(x) \leq \norm{f}_{\mathbb{X}} \norm{g}_{\mathbb{X^{'}}},
\end{equation}
for every $f \in \mathbb{X}$ and $g \in \mathbb{X}^{'}$. Every Banach function space 
coincides with its second associated space, that is, $\mathbb{X}=\mathbb{X}''$ (cf. \cite[Theorem 2.7, p.10]{BenSha88}). Consequently, 
\begin{equation}
\norm{f}_{ \mathbb{X}} = \sup \left\{ \int_{\Sigma} |f(x)g(x)| d\mu(x) 
: g \in  \mathbb{X'}, \enskip \norm{g}_{\mathbb{X'}} \leq 1 \right\}.
\end{equation}

The following theorem is proved in  \cite[Propositions 3.5 and 3.6, pp.\,15--16]{BenSha88}.

\begin{proposition}\label{prop:ACN-equiv}
Let $(\Sigma,\mu)$ be a $\sigma$-finite measure space and
$\mathbb{X}$ a Banach function space over $(\Sigma,\mu)$. The following statements are equivalent:
\begin{enumerate}
\item $\mathbb{X}$ has order continuous norm.
\item $\mathbb{X}$ has absolutely continuous norm: for every $f\in \mathbb{X}$ and every sequence of sets $\{E_n\}_{n=1}^{\infty} \subseteq \Sigma$ satisfying ${\bf 1}_{E_n}\to 0$ $\mu$-a.e. as $n\to\infty$ one has $\|f {\bf 1}_{E_n}\|_X\to 0$ as $n\to\infty$.
\item For every $f\in\mathbb{X}$, $g\in\mathbb{M}_\mu$, and every sequence $\{ f_n \}_n \subseteq \mathbb{M}_\mu$ satisfying 
$ |f_n|\leq |f|$ for every $n\geq 1$ and $f_n\to g \,\, \mu$-a.e. as  $n\to\infty$
one has that $\|f_n - g \|_{\mathbb{X}}\to 0$  as $n\to\infty$.
\end{enumerate}
\end{proposition}

We close this section by presenting an extrapolation theorem for the spaces $\mathbb{X}_w$ when $\mathbb{X}$ is a Banach function space. This theorem is a particular case of \cite[Theorem 3.1]{CaoMarMar22} and provides a generalization of the celebrated extrapolation theorem by Rubio de Francia.
Recall that a family of extrapolation pairs $\mathcal{F}$ is just
a collection of pairs of nonnegative measurable functions $(f,g)$.

\begin{theorem} \label{thm:extrapolation-BFS}
Let $(\Sigma,d,\mu)$ be a non-atomic space of homogeneous type
and let $\mathcal{F}$ be a family of extrapolation pairs. 
Fix two weights $v,w$ on $(\Sigma,\mu)$ such that $v(B)<\infty$ for every ball $B\subseteq\Sigma$.
Define $\mathcal{M}_v'h := \mathcal{M}(hv)/v$ for each $h\in\mathbb{M}_\mu$, where $\mathcal{M}$ is the Hardy-Littlewood maximal operator from \eqref{eq:maximal}.
Assume that $\mathbb{X}_v$ is a Banach function space over $(\Sigma,v)$ such that
\begin{align}
\label{eq:extrapolation-BFS-hyp1} \|(\mathcal{M}h)w\|_{\mathbb{X}_v} &\leq C_1 \|hw\|_{\mathbb{X}_v} \qquad \forall\,h \in \mathbb{M}_{\mu}, 
\\ \label{eq:extrapolation-BFS-hyp2} \|(\mathcal{M}'h)w^{-1} \|_{\mathbb{X}'_v} &\leq C_2 \|hw^{-1}\|_{\mathbb{X}'_v} \qquad \forall\,h \in\mathbb{M}_{\mu}.
\end{align}
Fix
$p_0\in(1,\infty)$ and suppose that there is a non-decreasing function $\Phi: [1, \infty) \to [1, \infty)$  such that for all $w_0\in A_{p_0}(\Sigma)$,
\begin{equation}  \label{eq:extrapolation-BFS-hyp3}
\|f\|_{L^{p_0}(\Sigma,w_0)} \leq \Phi([w_0]_{A_{p_0}})\|g\|_{L^{p_0}(\Sigma,w_0)}, \qquad\forall\,(f,g) \in \mathcal{F}.
\end{equation} 
Then,
\begin{equation}
\|fw\|_{\mathbb{X}_v} \leq C_0\, \|gw\|_{\mathbb{X}_v}, \qquad \forall\,(f,g) \in \mathcal{F},
\end{equation}
with $C_0=2^{3+2/p'_0}\,\Phi(2^{p_0}\,C_1^{p_0-1}\, C_2)$.
\end{theorem}

\section{The bounded approximation property in Banach function lattices} \label{sec:BAP}

Recall that a bounded linear operator $T:X \to Y$ between two Banach spaces $X$ and $Y$ is compact if it 
maps the unit ball of $X$ into a relatively compact set, that is, a set with compact closure in $Y$.
For example, one can easily check 
that finite rank operators (i.e., those with finite dimensional range) are compact. Moreover, 
the composition of a compact operator with a bounded operator is compact
and the norm-limit of a sequence of compact operators is compact. Consequently, if $\mathcal{K}(X\to Y)$ 
denotes the space of all bounded linear operators $K:X \to Y$ which are compact, then $\mathcal{K}(X\to Y)$ is closed
in the space of bounded linear operators $T:X \to Y$. Hence,
\begin{equation}
\|T\|^{\ess}_{X \to Y} = 0 \quad \Leftrightarrow \quad T \in \mathcal{K}(X\to Y),
\end{equation}
where $\|T\|^{\ess}_{X \to Y}$ denotes the essential norm, defined in \eqref{eq:ess-def}.

In this section, we identify conditions under which the identity operator can be uniformly approximated on compact sets by finite rank bounded linear operators, with uniformly bounded norms. This condition is made precise in the following definition.

\begin{definition}\label{def:BAP}
We say that a Banach space $X$ has the {\tt bounded approximation property} with constant $\lambda$ if 
for every compact set $F\subseteq X$ and every $\varepsilon>0$
there exists a finite rank bounded linear operator $P:X\to X$ such that
\begin{equation}
\sup_{f \in F} \|f - Pf\|_X < \varepsilon \quad \text{ and }
\quad \|P\|_{X\to X}\leq\lambda.
\end{equation}
\end{definition}

Recall the dyadic decomposition described in Theorem \ref{thm:dyadic-dec} for spaces of homogeneous type.
Choose an arbitrary 
$x_0 \in \Sigma$ and for every $n \in \mathbb{N}$, define 
\begin{equation} 
\mathbb{D}_n := \{Q^{n}_{\alpha} \in \mathbb{D}(\Sigma) :
\alpha\in I_n, \,\, d(Q^n_\alpha, x_0) \leq n \}. 
\end{equation} 
For every locally integrable function $f$, consider the operator 
\begin{equation} \label{eq:DefPN}
P_n f := \sum_{Q \in \mathbb{D}_n} \mathbf{1}_Q \fint_{Q} f d\mu. 
\end{equation}
Note that, as a consequence of Lebesgue's Differentiation Theorem for dyadic grids (cf. \cite[Corollary 2.6]{AndCruMoe15}),
\begin{equation}\label{eq:Pnf-to-f}
P_n f \to f \,\, \mu\text{-a.e. as } n \to \infty.
\end{equation}

\begin{proposition} \label{prop:PN}
Let $(\Sigma,d,\mu)$ be a space of homogeneous type, assume that $\mathbb{X}$ is a Banach function lattice over $(\Sigma,\mu)$ such that the Hardy-Littlewood maximal operator $\mathcal{M}$ is bounded on $\mathbb{X}$.
Then, for every $n\geq 1$, the operator $P_n:\mathbb{X}\to\mathbb{X}$ defined as in \eqref{eq:DefPN} is a finite rank bounded linear operator and
\begin{equation}\label{eq:PN-1}
\|P_n\|_{\mathbb{X}\to \mathbb{X}} \leq C_\Sigma \|\mathcal{M}\|_{\mathbb{X}\to \mathbb{X}} \qquad  \text{ for every } n \geq 1,
\end{equation}
where $C_\Sigma\in[1,\infty)$ is the constant from \eqref{eq:maximal-compar}.
Moreover, if $\mathbb{X}$ has order continuous norm, for every $f\in \mathbb{X}$ one has
\begin{equation} \label{eq:PNtoI}
\lim_{n \to \infty} \norm{f - P_n f}_{\mathbb{X}} = 0.
\end{equation}
\end{proposition}

\begin{proof}
Note that $|P_n f(x)| \leq \mathcal{M}_{\mathbb{D}} f(x)$ for every $x\in\Sigma$ and $n\geq 1$.
From this, \eqref{eq:maximal-compar}, and \eqref{eq:BL} we have
\begin{equation}
\|  P_n f  \|_{\mathbb{X}} \leq \|  \mathcal{M}_{\mathbb{D}}f  \|_{\mathbb{X}} \leq C_\Sigma
\|  \mathcal{M} f  \|_{\mathbb{X}} \leq C_\Sigma \|\mathcal{M}\|_{\mathbb{X}\to \mathbb{X}}
\| f  \|_{\mathbb{X}}, \qquad \forall f \in \mathbb{X}.
\end{equation}

Next, assume that $\mathbb{X}$ has order continuous norm and fix $f\in \mathbb{X}$. 
For every $n\geq 1$, define the sequence
\begin{equation}
g_n(x) := \sup_{m \geq n} |P_m f(x)-f(x)|.
\end{equation}
As $P_n f \to f$ $\mu$-a.e. (cf. \eqref{eq:Pnf-to-f}), we have that $g_n \to 0$ $\mu$-a.e. as $n\to\infty$.
Moreover, $g_n$ is decreasing and
\begin{equation}
g_n \leq \sup_{m\geq n}|P_m f| + |f|
\leq \mathcal{M}_{\mathbb{D}}f + |f| \leq (C_\Sigma+1) \mathcal{M}f \in \mathbb{X} \quad \mu\text{-a.e. on }\Sigma.
\end{equation}
As we are assuming that $\mathbb{X}$ has order continuous norm,
\begin{equation}
\|f - P_n f\|_{\mathbb{X}} \leq \|g_n\|_{\mathbb{X}} \to 0 \,\, \text{ as } \,\, n\to \infty.
\end{equation}
This completes the proof.
\end{proof}

The next theorem is the main result in this section, establishing the bounded approximation property in Banach function lattices under appropriate assumptions.

\begin{theorem}\label{thm:BAP-Xw}
Let $(\Sigma,d,\mu)$ be a space of homogeneous type and assume that $\mathbb{X}$ is a Banach function lattice over $(\Sigma,\mu)$ with order continuous norm such that the Hardy-Littlewood maximal operator $\mathcal{M}$ is bounded on $\mathbb{X}$. Consider the operators $P_n:\mathbb{X}\to\mathbb{X}$ defined in \eqref{eq:DefPN}. Then, for every compact set $F\subseteq \mathbb{X}$ and every $\varepsilon>0$ there is $n\geq 1$ such that
\begin{equation}
\sup_{f \in F} \|f - P_nf\|_{\mathbb{X}} < \varepsilon.
\end{equation}
As a consequence, $\mathbb{X}$ has the bounded approximation property with constant $C_\Sigma \|\mathcal{M}\|_{\mathbb{X}\to \mathbb{X}}$, where $C_\Sigma\in[1,\infty)$ is the constant from \eqref{eq:maximal-compar}.
\end{theorem}
\begin{proof}
Set $C_1 := C_\Sigma \|\mathcal{M}\|_{\mathbb{X}\to \mathbb{X}}\in[1,\infty)$.
Fix an arbitrary compact set $F\subseteq \mathbb{X}$ and $\varepsilon>0$.
As $F$ is compact, we can find 
functions $f_1,...,f_m \in \mathbb{X}$ such that 
\begin{equation} \label{eq:Ocovering-Xw}
F \subseteq \bigcup_{j=1}^{m} 
B_{\mathbb{X}}\Big(f_j,\frac{\varepsilon}{3C_1}\Big),
\end{equation}
where $B_{\mathbb{X}}$ denotes a ball in $\mathbb{X}$.
Using \eqref{eq:PNtoI}, there exists $n\geq 1$ large enough so that 
\begin{equation} \label{eq:epsilon00}
\norm{f_j - P_n f_j}_{\mathbb{X}}<\frac{\varepsilon}{3}
\qquad \text{ for every } j \in \{1,\dots,m\}
\end{equation}
For every arbitrary $f \in F$, pick $j_0\in\{1,\dots,m\}$ such 
that $\norm{f - f_{j_0}}_{\mathbb{X}}<\varepsilon/(3C_1)$. From this, and \eqref{eq:PN-1},
\begin{align}
\norm{f - P_n f}_{\mathbb{X}} 
&\leq \norm{f - f_{j_0}}_{\mathbb{X}} + \norm{f_{j_0} - P_n f_{j_0}}_{\mathbb{X}} 
+ \norm{P_n (f_{j_0} - f)}_{\mathbb{X}} 
\nonumber \\[4pt] 
& < \frac{\varepsilon}{3C_1} + \frac{\varepsilon}{3} 
+ C_1 \frac{\varepsilon}{3C_1}
\leq \varepsilon.
\end{align}
This proves that $\mathbb{X}$ has the bounded approximation property with constant $C_1$.
\end{proof}

The next result will be instrumental in the next section.

\begin{corollary}\label{cor:PNT}
Let $(\Sigma,d,\mu)$ be a space of homogeneous type and assume that $\mathbb{Y}$ is a Banach function lattice over $(\Sigma,\mu)$ with order continuous norm such that the Hardy-Littlewood maximal operator $\mathcal{M}$ is bounded on $\mathbb{Y}$. Let $X$ be an arbitrary Banach space and $T:X\to \mathbb{Y}$ a bounded linear operator. Consider the operators $P_n:\mathbb{Y}\to\mathbb{Y}$ defined in \eqref{eq:DefPN}. Then,
for every $\varepsilon>0$ there is $n\geq 1$ such that
\begin{equation}\label{eq:PN-d}
\| T-P_n T \| _{X\to\mathbb{Y}} < (2+C_\Sigma\|\mathcal{M}\|_{\mathbb{Y}\to \mathbb{Y}} )(\|T\|^{\ess}_{X\to\mathbb{Y}} + \varepsilon),
\end{equation}
where $C_\Sigma\in[1,\infty)$ is the constant from \eqref{eq:maximal-compar}.
\end{corollary}
\begin{proof}
Fix $\varepsilon>0$ arbitrarily and call 
\begin{equation}
\varepsilon_0 := \|T\|^{\ess}_{X\to\mathbb{Y}} + \varepsilon.
\end{equation}
Then, there exists a compact operator $K:X\to\mathbb{Y}$ such that 
$\norm{T-K}_{X\to\mathbb{Y}} < \varepsilon_0$. Define the set
$\mathcal{O} := K(B_{X}(0,1))$, where 
$B_{X}(0,1)$ denotes the unit ball in $X$. As $K$ is a compact operator,
$\overline{\mathcal{O}}$ is compact in $\mathbb{Y}$.
We now use Theorem \ref{thm:BAP-Xw}, which guarantees that there is $n\geq 1$ such that
$\|f-P_n f\|_{\mathbb{Y}}<\varepsilon_0$ for every $f\in\overline{\mathcal{O}}$. Therefore, we have $\norm{K - P_n K}_{X\to \mathbb{Y}}<\varepsilon_0$. As a consequence,
\begin{align}\label{eq:inproof1-PN}
\norm{T - P_n T}_{X\to \mathbb{Y}} 
& \leq \norm{T - K}_{X\to \mathbb{Y}} + \norm{K - P_n K}_{X\to \mathbb{Y}}
+ \norm{P_n (K - T)}_{X\to \mathbb{Y}} 
\nonumber \\[4pt] 
& < (2+C_\Sigma\|\mathcal{M}\|_{\mathbb{Y}\to \mathbb{Y}}) \varepsilon_0.
\end{align}
This finishes the proof.
\end{proof}

Note that in the previous result, $X$ is an arbitrary Banach space. In particular, the corollary is valid whenever $X=\mathbb{X}$ is a Banach function lattice. 
On the other hand, the spaces $\mathbb{X}_w$ introduced earlier inherit the Banach function lattice structure and the order continuity of the norm from $\mathbb{X}$ (cf. Proposition \ref{prop:Xw-properties}). This is particularly advantageous for the application of the results in this section to the spaces $\mathbb{X}_w$.

\section{Quantitative interpolation of compact operators} \label{sec:interpol}

\subsection{A general interpolation theorem for the essential norm}

%
We start with some basic definitions from interpolation theory.

\begin{definition}
A couple of Banach spaces $(X_0,X_1)$ is called {\tt compatible} if $X_0$ and $X_1$ are continuously embedded in the same Hausdorff topological vector space. In this context, a Banach space $X$ is {\tt intermediate} between $X_0$ and $X_1$ if 
\begin{equation}
X_0 \cap X_1 \subseteq X \subseteq X_0 + X_1,
\end{equation}
where the two inclusion maps are continuous.
\end{definition}

\begin{definition}
Let $(X_0,X_1)$ and $(Y_0,Y_1)$ be two compatible couples of Banach spaces
and $\Phi:[0,\infty) \times [0,\infty) \to [0,\infty)$ a function that is non-decreasing in each variable.
A pair of Banach spaces $(X,Y)$ is called a $\Phi$-{\tt interpolation} {\tt pair} between $(X_0,X_1)$ and $(Y_0,Y_1)$ if the following conditions are satisfied:
\begin{enumerate}
\item $X$ is intermediate between $X_0$ and $X_1$.
\item $Y$ is intermediate between $Y_0$ and $Y_1$.
\item If $T:X_0+X_1\to Y_0+Y_1$ is a linear operator and the restrictions $T:X_0\to Y_0$ and $T:X_1\to Y_1$ are both bounded, then $T:X\to Y$ is also bounded and
\begin{equation}
\|T\|_{X\to Y} \leq \Phi \left(\|T\|_{X_0\to Y_0}, \, \|T\|_{X_1 \to Y_1}\right).
\end{equation}
\end{enumerate}
\end{definition}

Of course, the main examples of interpolation pairs are given by the real and complex interpolation methods, for which $\Phi(x,y)=x^{1-\theta}y^{\theta}$ with $\theta\in(0,1)$.
We now turn to the proof of the general template for quantitative interpolation of compact operators announced in the introduction.

\begin{theorem}\label{thm:general}
Let $(\Sigma,d,\mu)$ be a space of homogeneous type 
and assume that $(X_0,X_1)$ and $(\mathbb{Y}_0,\mathbb{Y}_1)$ are two compatible couples of Banach spaces. 
Let $\Phi:[0,\infty) \times [0,\infty) \to [0,\infty)$ be a continuous function that is non-decreasing in each variable and
$(X,\mathbb{Y})$ a $\Phi$-interpolation pair between $(X_0,X_1)$ and $(\mathbb{Y}_0,\mathbb{Y}_1)$.
Suppose that $\mathbb{Y}_0$, $\mathbb{Y}_1$, and $\mathbb{Y}$ are Banach function lattices over $(\Sigma,\mu)$ and the Hardy-Littlewood maximal operator $\mathcal{M}$ is bounded on each of them.
Moreover, assume that $\mathbb{Y}_0$ and $\mathbb{Y}_1$ have order continuous norm.
Finally, fix a linear operator $T:X_0+X_1\to \mathbb{Y}_0+\mathbb{Y}_1$ such that the restrictions $T:X_0\to \mathbb{Y}_0$ and $T:X_1\to \mathbb{Y}_1$ are both bounded. Then,
\begin{equation}
\|T\|^{\ess}_{X\to \mathbb{Y}} \leq \Phi \left(C_0 \|T\|^{\ess}_{X_0\to \mathbb{Y}_0}, \, C_1 \|T\|^{\ess}_{X_1 \to \mathbb{Y}_1}\right),
\end{equation}
where 
\begin{equation}\label{eq:C0C1}
C_0 := 2+C_\Sigma \| \mathcal{M} \|_{\mathbb{Y}_0 \to \mathbb{Y}_0} \quad \text{ and } \quad C_1 := 2+C_\Sigma \| \mathcal{M} \|_{\mathbb{Y}_1 \to \mathbb{Y}_1},
\end{equation}
with $C_\Sigma$ the constant from \eqref{eq:maximal-compar}.

\end{theorem}
\begin{proof}
Since $(X,\mathbb{Y})$ is a $\Phi$-interpolation pair between $(X_0,X_1)$ and $(\mathbb{Y}_0,\mathbb{Y}_1)$, we have that $T:X\to \mathbb{Y}$ is bounded.
Moreover, as $\mathbb{Y}$ is a Banach function lattice where $\mathcal{M}$ is bounded, Proposition \ref{prop:PN} guarantees that $P_n:\mathbb{Y}\to \mathbb{Y}$ is a finite rank bounded linear operator.
As a consequence, $P_n T:X\to \mathbb{Y}$ is compact for every $n\geq 1$.
From this and Corollary \ref{cor:PNT}, for every $\varepsilon>0$ there is $n \geq 1$ such that
\begin{align}
\| T \|&^{\ess}_{X\to \mathbb{Y}}
 \leq \norm{T-P_n T}_{X\to \mathbb{Y}} 
\nonumber \\[4pt] 
&\qquad\leq \Phi( \norm{T-P_n T}_{X_0\to \mathbb{Y}_0} , \,
\norm{T-P_n T}_{X_1 \to \mathbb{Y}_1} )
\nonumber \\[4pt] 
&\qquad\leq \Phi( C_0 (\|T\|^{\ess}_{X_0\to\mathbb{Y}_0} + \varepsilon) , \,
C_1 (\|T\|^{\ess}_{X_1\to\mathbb{Y}_1} + \varepsilon) )
\end{align} 
with $C_0$ and $C_1$ as in \eqref{eq:C0C1}.
As we are assuming that $\Phi$ is continuous,
letting $\varepsilon \to 0$ we obtain the result.
\end{proof}

\subsection{Weighted Lebesgue spaces}

Consider a space of homogeneous type $(\Sigma,d,\mu)$ and a weight $w$ on $(\Sigma,\mu)$.
In this section we will apply the results in the previous sections to weighted Lebesgue spaces $L^p(\Sigma,w)$ with $p\in(1,\infty)$, defined in Section \ref{sec:prelim}.
From the Monotone Convergence Theorem we have that
\begin{equation}\label{eq:Lpw-hyp}
\parbox{9cm}{\begin{center}
$L^p(\Sigma,w)$ is a Banach function lattice over $(\Sigma,\mu)$ with order continuous norm.
\end{center}}
\end{equation}
Along with \eqref{eq:Mbounded-Lpw}, this guarantees that the results in Section \ref{sec:BAP} can be applied whenever $p\in(1,\infty)$ and $w\in A_p(\Sigma)$. In particular, $L^p(\Sigma,w)$ has the bounded approximation property in this setting.

Next, let $(\Sigma,\mu)$ be a $\sigma$-finite measure space, $p_0,p_1\in[1,\infty)$, and assume that $w_0, w_1$ are weights on $(\Sigma,\mu)$. Then,
\begin{equation}
L^{p_0}(\Sigma,w_0) \cap L^{p_1}(\Sigma,w_1)
\subseteq L^{p_\theta}(\Sigma,w_\theta) \subseteq
L^{p_0}(\Sigma,w_0) + L^{p_1}(\Sigma,w_1)
\end{equation}
with continuous inclusions, where
\begin{equation} 
w_\theta=w_0^{(1-\theta)p_\theta / p_0} w_1^{\theta p_\theta / p_1}, \qquad
\frac1{p_{\theta}} = \frac{1-\theta}{p_0} + \frac{\theta}{p_1}
\end{equation}
for every $\theta\in(0,1)$ (cf. \cite[Lemma 3.5, pp.211--212]{BenSha88}). This leads to the interpolation theorem with change of measures of Stein-Weiss. 
For the proof, the reader is referred to \cite[Theorem 3.6, p. 213]{BenSha88} (see also \cite{Ste56, SteWei58}).

\begin{theorem}\label{thm:stein-interpolation}
Let $(\Sigma,\mu)$ and $(\Sigma',\mu')$ be $\sigma$-finite measure spaces and assume that $w_0, w_1$ are weights on $(\Sigma,\mu)$ and $v_0, v_1$ are weights on $(\Sigma',\mu')$. Fix $p_0,p_1,q_0,q_1 \in [1,\infty)$
and let 
\begin{equation}
T:L^{p_0}(\Sigma,w_0)+L^{p_1}(\Sigma,w_1)\to L^{q_0}(\Sigma',v_0)+L^{q_1}(\Sigma',v_1)
\end{equation}
be a linear operator such that
\begin{align}
\norm{Tf}_{L^{q_0}(\Sigma',v_0)} \leq C_0 \norm{f}_{L^{p_0}(\Sigma,w_0)} 
\quad \forall f \in L^{p_0}(\Sigma,w_0), 
\\[4pt]
\norm{Tf}_{L^{q_1}(\Sigma',v_1)} \leq C_1 \norm{f}_{L^{p_1}(\Sigma,w_1)} 
\quad \forall f \in L^{p_1}(\Sigma,w_1),
\end{align}
for some $C_0, C_1 \in (0,\infty)$. Then, for every $\theta\in(0,1)$,
\begin{equation}
\norm{Tf}_{L^{q_{\theta}}(\Sigma',v_\theta)} 
\leq C_0^{1-\theta} C_1^{\theta} \norm{f}_{L^{p_{\theta}}(\Sigma,w_\theta)} 
\quad \forall f \in L^{p_{\theta}}(\Sigma,w_\theta),
\end{equation}
where
\begin{equation} 
w_\theta=w_0^{(1-\theta)p_\theta / p_0} w_1^{\theta p_\theta / p_1}, \qquad
\frac1{p_{\theta}} = \frac{1-\theta}{p_0} + \frac{\theta}{p_1},
\end{equation}
and
\begin{equation} \label{eq:vinterpol}
v_\theta=v_0^{(1-\theta)q_\theta / q_0} v_1^{\theta q_\theta / q_1}, \qquad
\frac1{q_{\theta}} = \frac{1-\theta}{q_0} + \frac{\theta}{q_1}.
\end{equation}
\end{theorem}

Note that, with the notation in \eqref{eq:vinterpol}, one may use Hölder's inequality to show that for every $\theta\in(0,1)$ and $q_0,q_1\in(1,\infty)$,
\begin{equation}
v_\theta\in A_{q_\theta}(\Sigma') \text{ whenever } 
v_0 \in A_{q_0}(\Sigma') \text{ and } v_1\in A_{q_1}(\Sigma').
\end{equation}

\begin{remark}\label{rem:RdF}
We can deduce Rubio de Francia extrapolation theorem from Theorem \ref{thm:extrapolation-BFS}. Indeed, if we set $v=1$, $\mathbb{X}_v=L^p(\Sigma,\mu)$ and apply the theorem with $w^{1/p}$ in place of $w$, then $\mathbb{X}'=L^{p'}(\Sigma,\mu)$ with $1/p+1/p'=1$ and \eqref{eq:extrapolation-BFS-hyp1}-\eqref{eq:extrapolation-BFS-hyp2} hold whenever $p\in(1,\infty)$ and $w\in A_p(\Sigma)$ (see \eqref{eq:Ap-p1} and \eqref{eq:Mbounded-Lpw}). In conclusion, if $(\Sigma, d,\mu)$ is a non-atomic space of homogeneous type and \eqref{eq:extrapolation-BFS-hyp3} holds, then for every $p\in(1,\infty)$ and $w\in A_p(\Sigma)$ one has
\begin{equation}
\|f\|_{L^{p}(\Sigma,w)} \leq C_0 \|g\|_{L^{p}(\Sigma,w)} \quad \text{ for every }(f,g) \in \mathcal{F}, 
\end{equation}
where $C_0=2^{3+2/p'_0}\,\Phi(C_{\Sigma,\mu,p_0,p}[w]^{1+(p_0-1)/(p-1)}_{A_{p}} )$.
\end{remark}

The following theorem is the main result of this section, addressing the quantitative interpolation of compact operators in the setting of weighted Lebesgue spaces.

\begin{theorem} \label{thm:compact-interpol}
Let $(\Sigma,d,\mu)$ and $(\Sigma',d',\mu')$ be spaces of homogeneous type  and fix four exponents
$p_0 , p_1\in [1,\infty)$ and $q_0,q_1 \in (1,\infty)$. Let $w_0$ and $w_1$ be arbitrary weights on $(\Sigma,\mu)$,
$v_0 \in A_{q_0}(\Sigma')$, and $v_1 \in A_{q_1}(\Sigma')$. 
Assume that
\begin{equation}
T:L^{p_0}(\Sigma,w_0)+L^{p_1}(\Sigma,w_1)\to L^{q_0}(\Sigma',v_0)+L^{q_1}(\Sigma',v_1)
\end{equation}
is a linear operator with the following bounded restrictions:
\begin{align}
T:L^{p_0}(\Sigma,w_0) \to L^{q_0}(\Sigma',v_0),
\\
T:L^{p_1}(\Sigma,w_1) \to L^{q_1}(\Sigma',v_1).
\end{align}
Then, for every $\theta\in(0,1)$ there is a constant $C\in[1,\infty)$ independent of T
such that
\begin{align} \label{eq:compact-interpol-1} 
\| T \|^{\ess}_{L^{p_{\theta}}(\Sigma,w_\theta)\to L^{q_{\theta}}(\Sigma',v_\theta)}
\leq C (\| T \|^{\ess}_{L^{p_0}(\Sigma,w_0)\to L^{q_0}(\Sigma',v_0)})^{1-\theta} 
\cdot (\| T \|^{\ess}_{L^{p_1}(\Sigma,w_1)\to L^{q_1}(\Sigma',v_1)})^{\theta},
\end{align}
where
\begin{equation} 
w_\theta=w_0^{(1-\theta)p_\theta / p_0} w_1^{\theta p_\theta / p_1}, \qquad
\frac1{p_{\theta}} = \frac{1-\theta}{p_0} + \frac{\theta}{p_1},
\end{equation}
and
\begin{equation} 
v_\theta=v_0^{(1-\theta)q_\theta / q_0} v_1^{\theta q_\theta / q_1}, \qquad
\frac1{q_{\theta}} = \frac{1-\theta}{q_0} + \frac{\theta}{q_1}.
\end{equation}
\end{theorem}

\begin{proof}
In view of \eqref{eq:Lpw-hyp}, \eqref{eq:Mbounded-Lpw}, and Theorem~\ref{thm:stein-interpolation}, we can apply Theorem \ref{thm:general} for weighted Lebesgue spaces with $\Phi(x,y)=x^{1-\theta}y^{\theta}$ to obtain
\begin{align}
\| T \|&^{\ess}_{L^{p_{\theta}}(\Sigma,w_\theta)\to L^{q_{\theta}}(\Sigma',v_\theta)}
 \leq (2+C_{\Sigma'} \| \mathcal{M} \|_{L^{q_0}(\Sigma',v_0)\to L^{q_0}(\Sigma',v_0)})^{1-\theta} \,(\| T \|^{\ess}_{L^{p_0}(\Sigma,w_0)\to L^{q_0}(\Sigma',v_0)})^{1-\theta}
\nonumber \\ & \qquad\quad\,\,\, \cdot (2+C_{\Sigma'} \| \mathcal{M} \|_{L^{q_1}(\Sigma',v_1)\to L^{q_1}(\Sigma',v_1)})^{\theta} \,(\| T \|^{\ess}_{L^{p_1}(\Sigma,w_1)\to L^{q_1}(\Sigma',v_1)})^{\theta}
\nonumber \\  & \qquad \leq \big(2+C_{\Sigma'} \, C_{\Sigma',\mu',q_0} \, [v_0]_{A_{q_0}}^{1/(q_0-1)} \big)^{1-\theta} \,(\| T \|^{\ess}_{L^{p_0}(\Sigma,w_0)\to L^{q_0}(\Sigma',v_0)})^{1-\theta}
\nonumber \\ & \qquad\quad\,\,\, \cdot \big(2+C_{\Sigma'} \, C_{\Sigma',\mu',q_1} \, [v_1]_{A_{q_1}}^{1/(q_1-1)} \big)^{\theta} \,(\| T \|^{\ess}_{L^{p_1}(\Sigma,w_1)\to L^{q_1}(\Sigma',v_1)})^{\theta},
\end{align} 
and this proves \eqref{eq:compact-interpol-1}.
\end{proof}

Building on the previous theorem, we prove two results that allow for the transfer of the control of the essential norm from unweighted to weighted Lebesgue spaces.

\begin{corollary}\label{cor:compact-interpol-1}
Let $(\Sigma,d,\mu)$ be a non-atomic space of homogeneous type and $p_0 \in (1,\infty)$.
Assume that $T$ is a linear operator such that for every weight $w_0\in A_{p_0}(\Sigma)$,
\begin{equation}
\|Tf\|_{L^{p_0}(\Sigma,w_0)} \leq \Phi([w_0]_{A_{p_0}})\|f\|_{L^{p_0}(\Sigma,w_0)}
\quad \forall f \in L^{p_0}(\Sigma,w_0), 
\end{equation}
where $\Phi: [1, \infty) \to [1, \infty)$ is a non-decreasing function.
Then, for every $p\in(1,\infty)$ and $w \in A_p(\Sigma)$ there is 
$\theta\in(0,1)$ such that $w^{1/\theta} \in A_p(\Sigma)$, and there is
$C\in[1,\infty)$ independent of T such that 
\begin{equation}
\|T\|^{\ess}_{L^p(\Sigma,w)\to L^p(\Sigma,w)} \leq C
(\|T\|^{\ess}_{L^{p}(\Sigma,\mu)\to L^{p}(\Sigma,\mu)})^{1-\theta}
(\|T\|^{\ess}_{L^{p}(\Sigma,w^{1/\theta})\to L^{p}(\Sigma,w^{1/\theta})})^{\theta}.
\end{equation}
As a consequence,
\begin{equation}
\|T\|^{\ess}_{L^p(\Sigma,w)\to L^p(\Sigma,w)} \leq C \,
\Phi( C_{\Sigma,\mu,p_0,p} [w^{1/\theta}]_{A_p}^{1+(p_0-1)/(p-1)})^\theta \,
(\|T\|^{\ess}_{L^{p}(\Sigma,\mu)\to L^{p}(\Sigma,\mu)})^{1-\theta}.
\end{equation}
\end{corollary}
\begin{proof}
It follows from Theorem \ref{thm:extrapolation-BFS} (see Remark \ref{rem:RdF})
that $T$ is bounded on $L^p(\Sigma,w)$ for all $p\in(1,\infty)$ and
$w \in A_p(\Sigma)$. Consequently, the result readily follows from \eqref{eq:openness} and
Theorem~\ref{thm:compact-interpol} with $p_0=p_1=q_0=q_1=p$, $v_0=w_0=1$, and $v_1=w_1=w^{1/\theta}$.
\end{proof}

\begin{corollary} \label{cor:compact-interpol-2}
Let $(\Sigma,d,\mu)$ be a non-atomic space of homogeneous type and fix two exponents $p_0,p_1 \in (1,\infty)$. Assume that $T$ is a linear operator such that for every weight $w_0 \in A_{p_0}(\Sigma)$, one has
\begin{equation} \label{eq:compact-interpol-2-hyp}
\|Tf\|_{L^{p_0}(\Sigma,w_0)} \leq \Phi([w_0]_{A_{p_0}})\|f\|_{L^{p_0}(\Sigma,w_0)}
\quad \forall f \in L^{p_0}(\Sigma,w_0), 
\end{equation}
where $\Phi: [1, \infty) \to [1, \infty)$ is a non-decreasing function.
Then, for every $p \in (1,\infty)$ and  $w \in A_p(\Sigma)$ there are $\theta\in(0,1)$ and $C\in[1,\infty)$ $($depending on $T$ only through $\Phi)$ 
such that
\begin{equation}
\|T\|^{\ess}_{L^{p}(\Sigma,w)\to L^{p}(\Sigma,w)} \leq C( \|T\|^{\ess}_{L^{p_1}(\Sigma,\mu)\to L^{p_1}(\Sigma,\mu)} )^{\theta}.
\end{equation}
In particular, if $T:L^{p_1}(\Sigma,\mu)\to L^{p_1}(\Sigma,\mu)$ is compact, then
\begin{equation}
T:L^p(\Sigma,w) \to L^p(\Sigma,w) \,\,\text{ is compact for every }
p \in (1,\infty) \text{ and } w \in A_p(\Sigma).
\end{equation}
\end{corollary}

\begin{proof}
Fix $p\in (1,\infty )$ and $w \in A_{p} (\Sigma )$. Pick $\theta _{1} \in (1/p, 1)$ sufficiently close to $1$ so that the exponent $q_{1}$ satisfying
\begin{equation}
\frac{1}{p} = \frac{1-\theta_1}{p_1} + \frac{\theta_1}{q_1}
\end{equation}
verifies $q_1\in(1,\infty)$.
From Corollary \ref{cor:compact-interpol-1} and Theorem~\ref{thm:compact-interpol}, there are $C\in[1,\infty)$ and $\theta_2\in(0,1)$ such that
$w^{1/\theta_2} \in A_p(\Sigma)$ and
\begin{align}
\|T\|^{\ess}_{L^p(\Sigma,w)\to L^p(\Sigma,w)} & \leq C
(\|T\|^{\ess}_{L^{p}(\Sigma,\mu)\to L^{p}(\Sigma,\mu)})^{1-\theta_2}
(\|T\|^{\ess}_{L^{p}(\Sigma,w^{1/\theta_2})\to L^{p}(\Sigma,w^{1/\theta_2})})^{\theta_2}
\nonumber \\[4pt] & \leq C
(\|T\|^{\ess}_{L^{p_1}(\Sigma,\mu)\to L^{p_1}(\Sigma,\mu)})^{(1-\theta_2)(1-\theta_1)}
(\|T\|^{\ess}_{L^{q_1}(\Sigma,\mu)\to L^{q_1}(\Sigma,\mu)})^{(1-\theta_2)\theta_1}
\nonumber \\[4pt] & \qquad \cdot 
(\|T\|^{\ess}_{L^{p}(\Sigma,w^{1/\theta_2})\to L^{p}(\Sigma,w^{1/\theta_2})})^{\theta_2}
\nonumber \\[4pt] & \leq C
(\|T\|^{\ess}_{L^{p_1}(\Sigma,\mu)\to L^{p_1}(\Sigma,\mu)})^{(1-\theta_2)(1-\theta_1)}
(\|T\|_{L^{q_1}(\Sigma,\mu)\to L^{q_1}(\Sigma,\mu)})^{(1-\theta_2)\theta_1}
\nonumber \\[4pt] & \qquad \cdot 
(\|T\|_{L^{p}(\Sigma,w^{1/\theta_2})\to L^{p}(\Sigma,w^{1/\theta_2})})^{\theta_2}.
\end{align}
Note that $T$ is bounded both on $L^{q_1}(\Sigma,\mu)$ and on $L^{p}(\Sigma,w^{1/\theta_2})$ as a consequence of \eqref{eq:compact-interpol-2-hyp} and  Theorem \ref{thm:extrapolation-BFS} (see also Remark \ref{rem:RdF}), and
\begin{align}
\|T\|^{\ess}_{L^p(\Sigma,w)\to L^p(\Sigma,w)} \leq C
C_1^{(1-\theta_2)\theta_1} C_2^{\theta_2}
(\|T\|^{\ess}_{L^{p_1}(\Sigma,\mu)\to L^{p_1}(\Sigma,\mu)})^{(1-\theta_2)(1-\theta_1)},
\end{align}
with
\begin{equation}
C_1=2^{3+2/p'_0}\,\Phi(C_{\Sigma,\mu,p_0,q_1}), \quad 
C_2=2^{3+2/p'_0}\,\Phi(C_{\Sigma,\mu,p_0,p}[w^{1/\theta_2}]^{1+(p_0-1)/(p-1)}_{A_{p}} ).
\end{equation}
Setting $\theta:=(1-\theta_2)(1-\theta_1)$ gives the result.
\end{proof}

\subsection{Rearrangement invariant Banach function spaces}\label{subsec:RIBFS}

One particularly relevant class of Banach function spaces is given by rearrangement invariant spaces.
Let $\mu_f$ denote the distribution function of a function $f \in \mathbb{M}_{\mu}$, that is, 
\begin{equation}
\mu_f (\lambda) = \mu ( \{ x \in \Sigma : |f(x)|>\lambda \} ) \text{,     } \qquad \lambda \geq 0. 
\end{equation}

\begin{definition}
A Banach function space $\mathbb{X}$ is called {\tt rearrangement} {\tt invariant} 
({\tt RIBFS}) if $\norm{ f } _{\mathbb{X}} = \norm{ g } _{\mathbb{X}}$ 
for every pair of functions $f,g \in \mathbb{X}$ verifying $\mu_f = \mu_g$.
\end{definition}

Many function spaces fall within this category, including Orlicz spaces (cf. Remark \ref{rem:RIBFS-examples} below).
Given a function $f \in \mathbb{M}_{\mu}$, the decreasing rearrangement of $f$ with 
respect to $\mu$ is the function $f^{*}_{\mu}$ defined as 
\begin{equation}
f^{*}_{\mu}(t) = \inf \{ \lambda \geq 0 : \mu_f (\lambda) \leq t \}, \quad t \in [0, \mu(\Sigma)). 
\end{equation} 
Note that $f$ and $f^{*}_{\mu}$ have the same distribution function. Consequently, one 
obtains
Luxemburg representation theorem: given a RIBFS $\mathbb{X}$, there exists another RIBFS
$\overline{\mathbb{X}}$ on $[0,\mu(\Sigma))$ such that $f \in \mathbb{X}$ if and only if 
$f_{\mu}^{*}\in \overline{\mathbb{X}}$ and 
$\norm{ f_{\mu}^{*} } _{\overline{\mathbb{X}}} = \norm{ f } _{\mathbb{X}}$. 
Given a weight $w$ on $(\Sigma,\mu)$,
these considerations naturally lead us to define a weighted RIBFS as:
\begin{equation} 
\mathbb{X}(w) := \{ f \in \mathbb{M}_{\mu} : \norm{ f^{*}_{w} } _{\overline{\mathbb{X}}} < \infty \}, 
\end{equation}
which is a Banach function space over $(\Sigma,w)$ with the norm 
$\norm{ f } _{\mathbb{X}(w)} := \norm{ f^{*}_{w} } _{\overline{\mathbb{X}}}$.
As an example, note that $\mathbb{X}=L^p(\Sigma,\mu)$ is a RIBFS over $(\Sigma,\mu)$ and $\mathbb{X}(w)=L^p(\Sigma,w)$.

Next, given a RIBFS $\mathbb{X}$, we define its lower and upper {\tt Boyd indices} as
\begin{equation}
p_{\mathbb{X}} := \lim_{t \to \infty} 
\frac{\log t}{\log \norm{ D_t } _{\overline{\mathbb{X}}\to \overline{\mathbb{X}}}} \text{,} \quad
\qquad q_{\mathbb{X}} 
:= \lim_{t \to 0^{+}} \frac{\log t}{\log \norm{ D_t } _{\overline{\mathbb{X}}\to \overline{\mathbb{X}}}}, 
\end{equation}
where $D_t: \overline{\mathbb{X}} \to \overline{\mathbb{X}}$ is the dilation operator defined by
\begin{equation}
D_t f(s):= \left\lbrace
\begin{array}{lcl}
f(s/t) & \text{ if } & s\in [0, \, \min \{ \mu(\Sigma), \, t \mu(\Sigma) \}),\\ 
0 & \text{ if } & s\in [ \min \{ \mu(\Sigma), \, t \mu(\Sigma) \} ,  \,\mu(\Sigma) ) .
\end{array}
\right. 
\end{equation}
Boyd indices capture essential information about the interpolation of operators in rearrangement invariant Banach function spaces (cf. \cite{BenSha88}).
The next result is a weighted version of Boyd's interpolation theorem which provides control 
over the operator norm and was proved in \cite[Theorem 8.24]{MarMarMitMitMit22}. 

\begin{theorem} \label{thm:Boyd}
Suppose that $\mathbb{X}$ is a RIBFS over a non-atomic, $\sigma$-finite measure space
$(\Sigma,\mu)$ with $\mu(\Sigma)>0$. Fix $p_0$ and $p_1$ so that
$1 \leq p_0 < p_{\mathbb{X}} \leq q_{\mathbb{X}} < p_1 \leq \infty$.
Let $w$ be a weight on $\Sigma$ such that $w(\Sigma)=\mu(\Sigma)$.
Suppose that
\begin{equation}
T:L^{p_0}(\Sigma,w)+L^{p_1}(\Sigma,w)\to L^{p_0}(\Sigma,w)+L^{p_1}(\Sigma,w)
\end{equation}
is a linear operator such that
\begin{equation} 
\norm{Tf}_{L^{p_0}(\Sigma,w)} \leq C_0 \norm{f}_{L^{p_0}(\Sigma,w)} \qquad \forall f \in L^{p_0}(\Sigma,w), 
\end{equation} 
\begin{equation} 
\norm{Tf}_{L^{p_1}(\Sigma,w)} \leq C_1 \norm{f}_{L^{p_1}(\Sigma,w)} \qquad \forall f \in L^{p_1}(\Sigma,w), 
\end{equation} 
for some $C_0, C_1 \in (0,\infty)$. Then,
\begin{equation}\label{eq:Xw-inclusion}
\mathbb{X}(w) \subseteq L^{p_0}(\Sigma,w) + L^{p_1}(\Sigma,w)
\end{equation}
and there exists $C\in(0,\infty)$ depending only on $p_0$ and $p_1$ such that
\begin{equation}
\norm{Tf}_{\mathbb{X}(w)} \leq C (C_0 + C_1) \norm{f}_{\mathbb{X}(w)} \qquad \forall f \in \mathbb{X}(w).
\end{equation} 
\end{theorem}

Boyd's interpolation theorem also holds for sublinear operators and one can use it to deduce (cf., e.g., \cite[(2.48)]{CaoMarMar22}) that if $\mathbb{X}$ is a RIBFS over a non-atomic space of homogeneous type
$(\Sigma,d,\mu)$, then
\begin{equation}\label{eq:M-bounded-Xw}
\parbox{13cm}{\begin{center}
$\mathcal{M}$ is bounded on $\mathbb{X}(w)$ whenever $1<p_{\mathbb{X}}\leq q_{\mathbb{X}} < \infty$ and $w\in A_{p_{\mathbb{X}}}(\Sigma)$.
\end{center} }
\end{equation}
Similarly, if we define $\mathcal{M}_v'h := \mathcal{M}(hv)/v$ for each $h\in\mathbb{M}_\mu$,
one can check (cf. \cite[(2.53)]{CaoMarMar22}) that
$\mathcal{M}'_v$ is bounded on $\mathbb{X}'(w)$ whenever $1<p_{\mathbb{X}}\leq q_{\mathbb{X}} < \infty$ and $w\in A_{p_{\mathbb{X}}}(\Sigma)$.
In turn, combining \eqref{eq:M-bounded-Xw} with Theorem \ref{thm:BAP-Xw} we deduce that
\begin{equation}
\parbox{12cm}{\begin{center}
$\mathbb{X}(w)$ has the bounded approximation property
whenever it has order continuous norm,
$1<p_\mathbb{X}\leq q_\mathbb{X} < \infty$, and 
$w\in A_{p_\mathbb{X}}(\Sigma)$.
\end{center} }
\end{equation}

\begin{remark}\label{rem:extrapolation-X}
Theorem \ref{thm:extrapolation-BFS} specializes to RIBFS by applying the result with $w$ in place of $v$, $1$ in place of $w$, and $\mathbb{X}(w)$ in place of $\mathbb{X}_v$.
As a consequence of the discussion above, hypotheses \eqref{eq:extrapolation-BFS-hyp1}-\eqref{eq:extrapolation-BFS-hyp2} hold whenever $1<p_{\mathbb{X}}\leq q_{\mathbb{X}} < \infty$ and $w\in A_{p_{\mathbb{X}}}(\Sigma)$.
We deduce that if $(\Sigma, d, \mu)$ is a non-atomic space of homogeneous type, $\mathbb{X}$ is a RIBFS over $(\Sigma,\mu)$ with $1<p_{\mathbb{X}}\leq q_{\mathbb{X}} < \infty$, $w\in A_{p_{\mathbb{X}}}(\Sigma)$, and \eqref{eq:extrapolation-BFS-hyp3} holds, then there is $C_0\in(0,\infty)$ such that
\begin{equation}
\|f\|_{\mathbb{X}(w)} \leq C_0 \, \|g\|_{\mathbb{X}(w)}, \qquad\forall\,(f,g) \in \mathcal{F}.
\end{equation}
\end{remark}

We are now prepared to prove our result on the interpolation of the essential norm on RIBFS.

\begin{theorem} \label{thm:compact-interpol-Xw}
Let $(\Sigma,d,\mu)$ be a non-atomic space of homogeneous type.
Let $p_0$ and $p_1$ be two exponents such that $1 < p_0 < p_1 < \infty$ and fix $w\in A_{p_0}(\Sigma)$ such that $w(\Sigma)=\mu(\Sigma)$. 
Let
\begin{equation}
T:L^{p_0}(\Sigma,w)+L^{p_1}(\Sigma,w)\to L^{p_0}(\Sigma,w)+L^{p_1}(\Sigma,w)
\end{equation}
be a linear operator with the following bounded restrictions:
\begin{align}
T:L^{p_0}(\Sigma,w) \to L^{p_0}(\Sigma,w),
\\
T:L^{p_1}(\Sigma,w) \to L^{p_1}(\Sigma,w).
\end{align}
Finally, assume that $\mathbb{X}$ is a RIBFS  over $(\Sigma,\mu)$
such that $p_0 < p_{\mathbb{X}} \leq q_{\mathbb{X}} < p_1$. Then, there exist $\theta,\theta'\in(0,1)$ and $C\in(0,\infty)$ independent of $T$ such that
\begin{align} \label{eq:compact-interpol-Xw-res}
\|T\|^{\ess}_{\mathbb{X}(w)\to \mathbb{X}(w)}
& \leq C \big( (\|T\|^{\ess}_{L^{p_0}(\Sigma,w)\to L^{p_0}(\Sigma,w)})^{1-\theta}
(\|T\|^{\ess}_{L^{p_1}(\Sigma,w)\to L^{p_1}(\Sigma,w)})^{\theta}
\nonumber \\ 
& \quad \,\,\,\, + (\|T\|^{\ess}_{L^{p_0}(\Sigma,w)\to L^{p_0}(\Sigma,w)})^{1-\theta'}
(\|T\|^{\ess}_{L^{p_1}(\Sigma,w)\to L^{p_1}(\Sigma,w)})^{\theta'} \big).
\end{align}
\end{theorem}

\begin{proof}
Note first that $T$ is bounded on $\mathbb{X}(w)$ as a consequence of Theorem \ref{thm:Boyd}.
As $w\in A_{p_0}(\Sigma) \subseteq A_{p_\mathbb{X}}(\Sigma)$, $\mathcal{M}$ is bounded on $\mathbb{X}(w)$ (see \eqref{eq:Ap-p2} and \eqref{eq:M-bounded-Xw}).
Next, choose $q_0, q_1$ so that $p_0<q_0<p_{\mathbb{X}}\leq q_{\mathbb{X}}<q_1<p_1$ and 
define $P_n:\mathbb{X}(w)\to\mathbb{X}(w)$ as in \eqref{eq:DefPN} for every $n\geq 1$. 
From Proposition \ref{prop:PN}, $P_n$ is a finite rank bounded linear operator on $\mathbb{X}(w)$ and hence $P_n T$ is compact in $\mathbb{X}(w)$.
Using Theorem~\ref{thm:Boyd} and Riesz-Thorin theorem, there are $C\in(0,\infty)$ and $\theta,\theta'\in(0,1)$ such that
\begin{align} \label{eq:compact-interpol-Xw-inproof1}
\|T\|^{\ess}_{\mathbb{X}(w)\to \mathbb{X}(w)} &\leq \norm{T-P_n T}_{\mathbb{X}(w)\to \mathbb{X}(w)} 
\nonumber \\[4pt] 
&\leq C \big( \norm{T-P_n T}_{L^{q_{0}}(\Sigma,w)\to L^{q_{0}}(\Sigma,w)} + \norm{T-P_n T} _{L^{q_{1}}(\Sigma,w)\to L^{q_{1}}(\Sigma,w)} \big)
\nonumber \\[4pt] 
&\leq C \big( \norm{T-P_n T}_{L^{p_{0}}(\Sigma,w)\to L^{p_{0}}(\Sigma,w)}^{1-\theta} \norm{T-P_n T}_{L^{p_{1}}(\Sigma,w)\to L^{p_{1}}(\Sigma,w)}^{\theta}  
\nonumber \\[4pt] 
& \quad + \norm{T-P_n T}_{L^{p_{0}}(\Sigma,w)\to L^{p_{0}}(\Sigma,w)}^{1-\theta'} 
\norm{T-P_n T}_{L^{p_{1}}(\Sigma,w)\to L^{p_{1}}(\Sigma,w)}^{\theta'} \big).
\end{align}
Finally, \eqref{eq:Lpw-hyp}, \eqref{eq:Ap-p2}, \eqref{eq:Mbounded-Lpw}, and Corollary \ref{cor:PNT} give the result.
\end{proof}

\begin{remark}
Naturally, a variant of Theorem \ref{thm:compact-interpol-Xw} can also be established with \eqref{eq:compact-interpol-Xw-res} replaced by
\begin{equation}
\|T\|^{\ess}_{\mathbb{X}(w)\to \mathbb{X}(w)}
\leq C (\|T\|^{\ess}_{L^{p_0}(\Sigma,w)\to L^{p_0}(\Sigma,w)} + 
\|T\|^{\ess}_{L^{p_1}(\Sigma,w)\to L^{p_1}(\Sigma,w)} ).
\end{equation}
However, this estimate does not allow one to obtain smallness for the essential norm on $\mathbb{X}(w)$ from smallness at a single Lebesgue endpoint; instead, it requires smallness at both endpoints.
\end{remark}

The following corollary is useful to transfer the compactness properties of a linear operator from (unweighted) Lebesgue spaces to weighted RIBFS.

\begin{corollary} \label{cor:compact-interpol-Xw-2}
Let $(\Sigma,d,\mu)$ be a non-atomic space of homogeneous type and fix two exponents $p_0,p_1\in(1,\infty)$. Assume that $T$ is a linear operator such that for every weight $w_0 \in A_{p_0}(\Sigma)$, one has
\begin{equation} \label{eq:compact-interpol-Xw-2-hyp}
\|Tf\|_{L^{p_0}(\Sigma,w_0)} \leq \Phi([w_0]_{A_{p_0}})\|f\|_{L^{p_0}(\Sigma,w_0)}
\quad \forall f \in L^{p_0}(\Sigma,w_0), 
\end{equation}
where $\Phi: [1, \infty) \to [1, \infty)$ is a non-decreasing function.
Let $\mathbb{X}$ be a RIBFS over $(\Sigma,\mu)$ such that 
$1<p_{\mathbb{X}} \leq q_{\mathbb{X}}<\infty$ and fix a weight $w \in A_{p_{\mathbb{X}}}(\Sigma)$ with $w(\Sigma)=\mu(\Sigma)$.
Then, there are $\theta,\theta'\in(0,1)$ and $C\in(0,\infty)$ $($depending on $T$ only through $\Phi)$ such that
\begin{equation}
\|T\|^{\ess}_{\mathbb{X}(w)\to \mathbb{X}(w)} \leq C \big( ( \|T\|^{\ess}_{L^{p_1}(\Sigma,\mu)\to L^{p_1}(\Sigma,\mu)} )^{\theta}
+( \|T\|^{\ess}_{L^{p_1}(\Sigma,\mu)\to L^{p_1}(\Sigma,\mu)} )^{\theta'} \big).
\end{equation}
In particular, if $T:L^{p_1}(\Sigma,\mu) \to L^{p_1}(\Sigma,\mu)$ is compact, then
\begin{equation}
T:\mathbb{X}(w) \longrightarrow \mathbb{X}(w) \quad \text{ is compact.} 
\end{equation}
\end{corollary}

\begin{proof}
Fix $w \in A_{p_{\mathbb{X}}}(\Sigma)$ with $w(\Sigma)=\mu(\Sigma)$. As a consequence of \eqref{eq:Ap-p2} and \eqref{eq:openness}, $w \in A_{p_{\mathbb{X}}-\varepsilon}(\Sigma)$ 
and $w \in A_{q_{\mathbb{X}}+\varepsilon}(\Sigma)$ for some $\varepsilon>0$ with 
$p_{\mathbb{X}}-\varepsilon \in(1,\infty)$. 
Next, we apply Corollary \ref{cor:compact-interpol-2} to obtain
$\theta_1,\theta_2\in(0,1)$ and $C_1,C_2\in[1,\infty)$ such that
\begin{equation}
\|T\|^{\ess}_{L^{p_{\mathbb{X}}-\varepsilon}(\Sigma,w)\to L^{p_{\mathbb{X}}-\varepsilon}(\Sigma,w)} \leq C_1( \|T\|^{\ess}_{L^{p_1}(\Sigma,\mu)\to L^{p_1}(\Sigma,\mu)} )^{\theta_1}
\end{equation}
and 
\begin{equation}
\|T\|^{\ess}_{L^{q_{\mathbb{X}}+\varepsilon}(\Sigma,w)\to L^{q_{\mathbb{X}}+\varepsilon}(\Sigma,w)} \leq C_2( \|T\|^{\ess}_{L^{p_1}(\Sigma,\mu)\to L^{p_1}(\Sigma,\mu)} )^{\theta_2}.
\end{equation}
Finally, applying  Theorem \ref{thm:compact-interpol-Xw}, that there exist $\theta_0,\theta_0'\in(0,1)$ and $C\in(0,\infty)$ such that
\begin{align}
\|T\|^{\ess}_{\mathbb{X}(w)\to \mathbb{X}(w)}
& \leq C \big( (\|T\|^{\ess}_{L^{p_{\mathbb{X}}-\varepsilon}(\Sigma,w)\to L^{p_{\mathbb{X}}-\varepsilon}(\Sigma,w)})^{1-\theta_0}
(\|T\|^{\ess}_{L^{q_{\mathbb{X}}+\varepsilon}(\Sigma,w)\to L^{q_{\mathbb{X}}+\varepsilon}(\Sigma,w)})^{\theta_0}
\nonumber\\ & \quad \,\, + 
(\|T\|^{\ess}_{L^{p_{\mathbb{X}}-\varepsilon}(\Sigma,w)\to L^{p_{\mathbb{X}}-\varepsilon}(\Sigma,w)})^{1-\theta'_0}
(\|T\|^{\ess}_{L^{q_{\mathbb{X}}+\varepsilon}(\Sigma,w)\to L^{q_{\mathbb{X}}+\varepsilon}(\Sigma,w)})^{\theta'_0} \big)
\nonumber\\ & \leq 
C \, C_1^{1-\theta_0} C_2^{\theta_0}  ( \|T\|^{\ess}_{L^{p_1}(\Sigma,\mu)\to L^{p_1}(\Sigma,\mu)} )^{\theta_1(1-\theta_0) + \theta_2\theta_0}.
\nonumber\\ & \quad \,\, + 
C \, C_1^{1-\theta'_0} C_2^{\theta'_0}  ( \|T\|^{\ess}_{L^{p_1}(\Sigma,\mu)\to L^{p_1}(\Sigma,\mu)} )^{\theta_1(1-\theta'_0) + \theta_2\theta'_0}.
\end{align}
Setting $\theta:=\theta_1(1-\theta_0) + \theta_2\theta_0$ and $\theta':=\theta_1(1-\theta'_0) + \theta_2\theta'_0$,  we obtain the result.
\end{proof}

\begin{remark} \label{rem:RIBFS-examples}
The results established in this section apply to a wide range of function spaces. Apart from ordinary and weighted Lebesgue spaces, it includes the cases
\begin{equation}
\mathbb{X}=L^p(\Sigma,\mu) + L^q(\Sigma,\mu)
\quad\text{ and }\quad\mathbb{X}=L^p(\Sigma,\mu) \cap L^q(\Sigma,\mu)
\end{equation}
with $p,q\in(1,\infty)$, as well as their weighted counterparts.
It also includes other Orlicz spaces such as Zygmund spaces 
\begin{equation}
\mathbb{X}=L^{p}(\log L)^{\alpha}(\Sigma,\mu),
\end{equation}
with $p\in(1,\infty)$ and $\alpha \in \mathbb{R}$, and their weighted versions
(cf., e.g., \cite[Chapter 4]{CruMarPer11}).
\end{remark}

\subsection{Variable exponent Lebesgue spaces}\label{subsec:Lpvar}

In this section, we consider variable exponent Lebesgue spaces. In this context, given a $\sigma$-finite measure space
$(\Sigma,\mu)$, measurable functions $p:\Sigma \to (1,\infty)$ are called {\tt exponent} {\tt functions}.
We will denote these functions by $p(\cdot)$ to highlight the fact that the exponent is variable. To describe their range, we define 
\begin{equation} 
\displaystyle p_{-} := \essinf_{x \in \Sigma} p(x), \qquad p_{+} := \esssup_{x \in \Sigma} p(x). 
\end{equation}
Given an exponent function $p(\cdot)$, the {\tt variable} {\tt exponent} {\tt Lebesgue} {\tt space}  $L^{p(\cdot)}(\Sigma,\mu)$ is the Banach space defined through the norm
\begin{equation}\label{eq:Lpvar-def}
\|f\|_{L^{p(\cdot)}(\Sigma,\mu)} 
:= \inf \left\lbrace \lambda >0: \,\, \int_{\Sigma}
\left( \frac{|f(x)|}{\lambda} \right)^{p(x)}
d\mu(x) \leq 1 \right\rbrace .
\end{equation}
Note that if $p(\cdot)=p$ is constant, we obtain the classical Lebesgue space $L^p(\Sigma,\mu)$. Next, we recall that weighted variable exponent Lebesgue spaces $L^{p(\cdot)}(\Sigma,w)$ are defined through the norm
\begin{equation}
\| f \| _{L^{p(\cdot)}(\Sigma,w)} :=
\| f w \| _{L^{p(\cdot)}(\Sigma,\mu)}.
\end{equation}

Variable exponent Lebesgue  spaces have been extensively studied for the Euclidean case (cf., e.g.,  \cite{CruFio13}). 
For the case of spaces of homogeneous type, the reader is referred to \cite{AdaHarHas15,DieHarHasRuz11,HarHasPer04}. 
The space $L^{p(\cdot)}(\Sigma,\mu)$ is a Banach function space (cf. \cite[Theorem 3.2.13]{DieHarHasRuz11}) but in general, it is not rearrangement invariant.
Moreover, whenever $\mu$ is $\sigma$-finite, non-atomic, and complete, one has that
\begin{equation}
\parbox{13cm}{\begin{center}
$L^{p(\cdot)}(\Sigma,\mu)$ has absolutely continuous norm if and only if $p_+<\infty$.
\end{center} }
\end{equation}
(cf. \cite[Lemma 3.1]{LukPikPok11}). Combining this with Proposition \ref{prop:ACN-equiv},
we deduce that $L^{p(\cdot)}(\Sigma,\mu)$ has order continuous norm if and only if $p_+<\infty$.
In turn, from this and Proposition \ref{prop:Xw-properties}, we obtain that for every weight $w$ on $(\Sigma,\mu)$,
\begin{equation}\label{eq:Lpvarw-abscont}
\parbox{8.5cm}{\begin{center}
$L^{p(\cdot)}(\Sigma,w)$ is a Banach function lattice with order continuous norm if and only if $p_+<\infty$.
\end{center} }
\end{equation}

Given a space of homogeneous type $(\Sigma,d,\mu)$, we say that an exponent function $p(\cdot)$ is {\tt log-H\"older} {\tt continuous} on $\Sigma$, denoted as $p(\cdot) \in {\rm LH}(\Sigma)$, if there are $C_1,C_2 \in (0,\infty)$, $p_{\infty} \in \mathbb{R}$, and a basepoint $x_0\in\Sigma$ such that
\begin{equation}
|p(x) - p(y)| < \frac{-C_1}{\log(d(x , y))}
\quad \text{ whenever } d(x,y)<\frac{1}{2}
\end{equation}
and
\begin{equation} \label{eq:LH-inf}
|p(x) - p_{\infty}| < \frac{C_2}{\log (e + d(x,x_0))}
\quad \text{ for every } x\in\Sigma.
\end{equation}

One can check that the choice of the basepoint $x_0\in\Sigma$ above is irrelevant, that is, if \eqref{eq:LH-inf} holds for some $x_0\in\Sigma$ then it also holds for any other $y_0\in\Sigma$, possibly with a different constant (cf. \cite[Lemma 2.1]{AdaHarHas15}). Moreover, $p(\cdot)\in{\rm LH}(\Sigma)$ is equivalent to $1/p(\cdot)\in{\rm LH}(\Sigma)$ whenever $p(\cdot)$ is bounded.

Next, we define the class of $A_{p(\cdot)}(\Sigma)$ weights on $(\Sigma,\mu)$ as the family of weights satisfying
\begin{equation}\label{eq:Apvar}
[w]_{A_{p(\cdot)}}:=\sup_{B} \big( [\mu(B)]^{-1}\,
\norm{\mathbf{1}_B w}_{L^{p(\cdot)}(\Sigma,\mu)}
\|\mathbf{1}_B w^{-1}\|_{L^{p'(\cdot)}(\Sigma,\mu)}\big)<\infty,
\end{equation}
where the supremum runs over all balls $B$ in $\Sigma$ and $1/p(\cdot)+1/p'(\cdot)=1$.
Under the background assumption that $1<p_{-}\leq p_{+}<\infty$ and $p(\cdot)\in{\rm LH}(\Sigma)$, it is proved in \cite[Theorem 1.13]{CruCum22} that
\begin{equation}\label{eq:Apvar-equiv}
w \in A_{p(\cdot)}(\Sigma) \,\, \Longleftrightarrow \,\, \|\mathcal{M}f\|_{L^{p(\cdot)}(\Sigma,w)} \leq C \|f\|_{L^{p(\cdot)}(\Sigma,w)}.
\end{equation}
In particular, $L^{p(\cdot)}(\Sigma,w)\subseteq L^1_{\rm loc}(\Sigma,\mu)$ whenever $1<p_{-}\leq p_{+}<\infty$,  $p(\cdot)\in{\rm LH}(\Sigma)$, and $w \in A_{p(\cdot)}(\Sigma)$.
On the other hand, if $p(\cdot)\in{\rm LH}(\Sigma)$ and $p_{-}>1$ then $1\in A_{p(\cdot)}$ (cf. \cite[Lemma 3.6]{CruCum22}), hence, we recover the corresponding result in the unweighted setting:
\begin{equation}\label{eq:M-Lpvar}
\|\mathcal{M}f\|_{L^{p(\cdot)}(\Sigma,\mu)} \leq C \|f\|_{L^{p(\cdot)}(\Sigma,\mu)} \quad \text{ if } 1<p_{-}\leq p_{+}<\infty \text{ and } p(\cdot)\in{\rm LH}(\Sigma).
\end{equation}
Similarly, $\|\mathcal{M}f\|_{L^{p'(\cdot)}(\Sigma,\mu)} \leq C \|f\|_{L^{p'(\cdot)}(\Sigma,\mu)}$ under the same hypotheses.

Combining the results above with Theorem \ref{thm:BAP-Xw} we immediately obtain the following theorem.

\begin{theorem}\label{thm:BAP-Lpvar}
Let $(\Sigma,d,\mu)$ be a non-atomic, complete space of homogeneous type.
Let $p(\cdot)\in {\rm LH}(\Sigma)$ be an exponent function such that $1<p_{-} \leq p_{+}<\infty$ and $w\in A_{p(\cdot)}(\Sigma)$.
Then, $L^{p(\cdot)}(\Sigma,w)$ has the bounded approximation property.
\end{theorem}

With regard to interpolation in variable Lebesgue spaces, the following result was established in \cite[Section 7.1]{DieHarHasRuz11} for the special case $p_i(\cdot)=q_i(\cdot)$. As noted in \cite[p. 128]{CruFio13}, the more general statement presented below can be derived by an identical proof.

\begin{theorem} \label{thm:Lpvar-interpolation}
Let $(\Sigma,\mu)$ be a complete, $\sigma$-finite measure space with $\mu(\Sigma)>0$ and fix four exponent functions
$p_0(\cdot)$, $p_1(\cdot)$, $q_0(\cdot)$, and $q_1(\cdot)$. Suppose that
\begin{equation}
T:L^{p_0(\cdot)}(\Sigma,\mu)+L^{p_1(\cdot)}(\Sigma,\mu)\to L^{q_0(\cdot)}(\Sigma,\mu) + L^{q_1(\cdot)}(\Sigma,\mu)
\end{equation}
is a linear operator such that 
\begin{equation}
\norm{Tf}_{L^{q_0(\cdot)}(\Sigma,\mu)} \leq C_0 \norm{f}_{L^{p_0(\cdot)}(\Sigma,\mu)} 
\quad \forall f \in L^{p_0(\cdot)}(\Sigma,\mu), 
\end{equation}
\begin{equation}
\norm{Tf}_{L^{q_1(\cdot)}(\Sigma,\mu)} \leq C_1 \norm{f}_{L^{p_1(\cdot)}(\Sigma,\mu)} 
\quad \forall f \in L^{p_1(\cdot)}(\Sigma,\mu), 
\end{equation}
for some $C_0,C_1\in(0,\infty)$. Then, for every $\theta\in(0,1)$ there is $C\in(0,\infty)$ such that
\begin{equation} 
\norm{Tf}_{L^{q_{\theta}(\cdot)}(\Sigma,\mu)} \leq C C_0 ^{1-\theta} C_1 ^{\theta} 
\norm{f}_{L^{p_{\theta}(\cdot)}(\Sigma,\mu)} 
\quad \forall f \in L^{p_\theta(\cdot)}(\Sigma,\mu), 
\end{equation}
where
\begin{equation} 
\frac1{p_{\theta}(x)} = \frac{1-\theta}{p_0(x)} + \frac{\theta}{p_1(x)}, \qquad
\frac1{q_{\theta}(x)} = \frac{1-\theta}{q_0(x)} + \frac{\theta}{q_1(x)}.
\end{equation}
\end{theorem}

\begin{remark}\label{rem:extrapolation-Lpvar}
Theorem \ref{thm:extrapolation-BFS} also applies in the setting of variable exponent Lebesgue spaces. Although the result is applicable to the weighted case, for our purposes it suffices to consider the unweighted setting.
To this end, we choose $v=w=1$ and $\mathbb{X}=L^{p(\cdot)}(\Sigma,\mu)$, so that $\mathbb{X}'=L^{p'(\cdot)}(\Sigma,\mu)$ with $1/p(\cdot)+1/p'(\cdot)=1$. Then, hypotheses \eqref{eq:extrapolation-BFS-hyp1}-\eqref{eq:extrapolation-BFS-hyp2} hold whenever $1<p_{-}\leq p_{+}<\infty$  and $p(\cdot)\in{\rm LH}(\Sigma)$ (see \eqref{eq:M-Lpvar}). 
We deduce that if $(\Sigma, d,\mu)$ is a non-atomic space of homogeneous type, $p(\cdot)\in{\rm LH}(\Sigma)$ verifies $1<p_{-}\leq p_{+}<\infty$, and \eqref{eq:extrapolation-BFS-hyp3} holds, then there is $C_0\in(0,\infty)$ such that
\begin{equation}
\|f\|_{L^{p(\cdot)}(\Sigma,\mu)} \leq C_0 \|g\|_{L^{p(\cdot)}(\Sigma,\mu)} \quad \text{ for every }(f,g) \in \mathcal{F}.
\end{equation}
\end{remark}

The next theorem presents the main result of this section, addressing quantitative interpolation for compact operators in variable exponent Lebesgue spaces.

\begin{theorem} \label{thm:compact-var-interpol}
Let $(\Sigma,d,\mu)$ be a non-atomic, complete space of homogeneous type. Assume that
$p_0(\cdot)$ and  $p_1(\cdot)$ are arbitrary exponent functions and fix $q_0(\cdot)\in{\rm LH}(\Sigma)$ and $q_1(\cdot)\in{\rm LH}(\Sigma)$ satisfying
\begin{equation}
1<(q_0)_{-} \leq (q_0)_{+} <\infty, \qquad 1<(q_1)_{-} \leq (q_1)_{+} <\infty.
\end{equation}
Let 
\begin{equation}
T:L^{p_0(\cdot)}(\Sigma,\mu)+L^{p_1(\cdot)}(\Sigma,\mu)\to L^{q_0(\cdot)}(\Sigma,\mu) + L^{q_1(\cdot)}(\Sigma,\mu)
\end{equation}
be a linear operator with the following bounded restrictions:
\begin{align}
T:L^{p_0(\cdot)}(\Sigma,\mu) \to L^{q_0(\cdot)}(\Sigma,\mu),
\\
T:L^{p_1(\cdot)}(\Sigma,\mu) \to L^{q_1(\cdot)}(\Sigma,\mu).
\end{align}
Then, for every $\theta\in(0,1)$ there is a constant
$C\in[1,\infty)$ independent of $T$ such that
\begin{align} \label{eq:compact-var-interpol-1} 
\| T \|^{\ess}_{L^{p_{\theta}(\cdot)}(\Sigma,\mu)\to L^{q_{\theta}(\cdot)}(\Sigma,\mu)}
\leq  C(\| T \|^{\ess}_{L^{p_0(\cdot)}(\Sigma,\mu)\to L^{q_0(\cdot)}(\Sigma,\mu)})^{1-\theta} 
\cdot (\| T \|^{\ess}_{L^{p_1(\cdot)}(\Sigma,\mu)\to L^{q_1(\cdot)}(\Sigma,\mu)})^{\theta},
\end{align}
where
\begin{equation} 
\frac1{p_{\theta}(x)} = \frac{1-\theta}{p_0(x)} + \frac{\theta}{p_1(x)}, \qquad
\frac1{q_{\theta}(x)} = \frac{1-\theta}{q_0(x)} + \frac{\theta}{q_1(x)}.
\end{equation}
\end{theorem}

\begin{proof}
It is a direct consequence of Theorem \ref{thm:general} with $\Phi(x,y)=x^{1-\theta}y^{\theta}$, whose applicability is ensured by \eqref{eq:Lpvarw-abscont}, \eqref{eq:M-Lpvar}, and Theorem \ref{thm:Lpvar-interpolation}.
\end{proof}

For subsequent applications, we formulate the following corollary, which enables the transfer of the smallness of the essential norm from classical Lebesgue spaces to variable exponent Lebesgue spaces.

\begin{corollary} \label{cor:compact-interpol-Lpvar-2}
Let $(\Sigma,d,\mu)$ be a non-atomic, complete space of homogeneous type and $p_1\in(1,\infty)$.
Assume that $T$ is a linear operator such that 
\begin{equation}
\parbox{12cm}{\begin{center}
$T:L^{p_0(\cdot)}(\Sigma,\mu)\to L^{p_0(\cdot)}(\Sigma,\mu)$ is bounded for every $p_0(\cdot) \in {\rm LH}(\Sigma)$ with $1<(p_0)_{-}\leq (p_0)_{+}<\infty$.
\end{center} }
\end{equation}
Then, for every exponent function $p(\cdot)\in {\rm LH}(\Sigma)$ satisfying $1< p_{-} \leq p_{+} < \infty$,
there are $q(\cdot)\in{\rm LH}(\Sigma)$ with $1<q_{-}\leq q_{+} < \infty$, $\theta\in(0,1)$ and $C\in(0,\infty)$ independent of $T$ such that 
\begin{align}
\| T \|^{\ess}_{L^{p(\cdot)}(\Sigma,\mu)\to L^{p(\cdot)}(\Sigma,\mu)}
\leq  C(\| T \|^{\ess}_{L^{p_1}(\Sigma,\mu)\to L^{p_1}(\Sigma,\mu)})^{1-\theta} 
\cdot (\| T \|^{\ess}_{L^{q(\cdot)}(\Sigma,\mu)\to L^{q(\cdot)}(\Sigma,\mu)})^{\theta}.
\end{align}
In particular, if
$T:L^{p_1}(\Sigma,\mu)\to L^{p_1}(\Sigma,\mu)$ is compact, then
\begin{equation}
\parbox{10cm}{\begin{center}
$T:L^{p(\cdot)}(\Sigma,\mu)\to L^{p(\cdot)}(\Sigma,\mu)$ is compact for every $p(\cdot) \in {\rm LH}(\Sigma)$ with $1<p_{-}\leq p_{+}<\infty$.
\end{center} }
\end{equation}
\end{corollary}
\begin{proof}
Fix an exponent function $p(\cdot) \in {\rm LH}(\Sigma)$ 
satisfying $1<p_{-}\leq p_{+}<\infty$. Choose $q(\cdot)$ so that
\begin{equation}
\frac1{p(x)} = \frac{1-\theta}{p_1} + \frac{\theta}{q(x)}
\end{equation}
for some $\theta\in(0,1)$.
Since $p(\cdot) \in {\rm LH}(\Sigma)$, one can check that $q(\cdot) \in {\rm LH}(\Sigma)$, 
and taking $\theta$ suitably close to $1$, we can also guarantee that $1<q_{-}\leq q_{+} < \infty$. Then, $T$ is bounded both on $L^{p_1}(\Sigma,\mu)$ and on $L^{q(\cdot)}(\Sigma,\mu)$ by hypothesis and the result follows at once from Theorem \ref{thm:compact-var-interpol}.
\end{proof}

\section{Applications}\label{sec:applications}

\subsection{Essential norm of Calderón-Zygmund commutators}

The aim of this section is to bound the essential norm of commutators $[M_b,T]f := bT(f) - T(bf)$ in the context of spaces of homogeneous type $(\Sigma,d,\mu)$, where $T$ is a Calderón-Zygmund operator and $b \in \mathrm{BMO}(\Sigma,\mu)$.

First, recall that the space of functions with bounded mean oscillations, denoted by $\BMO(\Sigma,\mu)$, is the space defined through the semi-norm 
\begin{equation}\label{eq:BMO-seminorm1}
\norm{f}_{\BMO(\Sigma,\mu)}
:=\sup_{B\subseteq\Sigma}\fint_{B}\Big| f 
-\fint_B f\,d\mu\Big|\,d\mu \qquad \text{ if } \mu(\Sigma)=+\infty,
\end{equation}
where the supremum is taken over all balls $B$ in $\Sigma$, and
\begin{equation}\label{eq:BMO-seminorm2}
\norm{f}_{\BMO(\Sigma,\mu)}
:=\sup_{B\subseteq\Sigma}\fint_{B}\Big| f 
-\fint_B f\,d\mu\Big|\,d\mu +\Big| \int_{\Sigma}f\,d\mu\Big|  \qquad \text{ if } \mu(\Sigma)<+\infty.
\end{equation}
The space of functions with vanishing mean oscillations $\VMO(\Sigma,\mu)$ is defined as the closure in $\BMO(\Sigma,\mu)$ of the space of uniformly continuous $\BMO(\Sigma,\mu)$ functions. If $(\Sigma,d,\mu)$ is compact then there is $\alpha_0\in(0,1)$ such that $\VMO(\Sigma,\mu)$ is the closure of the Hölder space $\mathscr{C}^{\alpha}(\Sigma,d)$ in $\BMO(\Sigma,\mu)$ for every $\alpha\in(0,\alpha_0)$ (cf. \cite[Proposition 2.15]{HofMitTay09}).

\begin{definition}
Let $(\Sigma,d,\mu)$ be a space of homogeneous type.
We say that a linear continuous operator
$T:\mathscr{C}^{\alpha}(\Sigma,d)\to (\mathscr{C}^{\alpha}(\Sigma,d))^{*}$ with $\alpha\in(0,\alpha_0)$ has an associated 
kernel $k \in L^1_{\rm loc}(\Sigma\times\Sigma\setminus{\rm diag}(\Sigma))$, 
where ${\rm diag}(\Sigma)=\{(x,x):x\in\Sigma\}$, if
\begin{equation}\label{eq:CZO-def-1}
\langle Tf , g \rangle = \int_\Sigma \int_\Sigma
k(x,y) f(y) g(x) \, d\mu(x) \, d\mu(y),
\end{equation}
whenever $f,g\in \mathscr{C}^{\alpha}(\Sigma,d)$ have bounded, disjoint supports.
Moreover, $T$ is called a {\tt Calder\'on-Zygmund} {\tt operator} if it 
is bounded on $L^2(\Sigma,\mu)$ and
has an associated kernel $k$
such that there are $C\in(0,\infty)$ and $\varepsilon>0$ with
\begin{equation}
|k(x,y)| \le \frac{C}{\mu\bigl(B(x, d(x,y))\bigr)},
\qquad \forall\, x,y \in \Sigma,
\end{equation}
and
\begin{equation}
|k(x,y) - k(x',y)|  + |k(y,x) - k(y,x')| 
\leq
\frac{C}{\mu\bigl(B(x, d(x,y))\bigr)}
\left(\frac{d(x,x')}{d(x,y)}\right)^{\varepsilon}
\end{equation}
for all $x, x', y \in \Sigma$ satisfying $d(x,x') \leq \varepsilon d(x,y)$.
\end{definition}

It is well-known that a Calderón-Zygmund operator, which is by definition bounded on $L^2(\Sigma,\mu)$, is also bounded on $L^p(\Sigma,w)$ for every $p\in(1,\infty)$ and $w\in A_p(\Sigma)$ (cf. \cite{Mar04} for a proof in the context of spaces of homogeneous type). Moreover, Coifman, Rochberg and Weiss proved that commutators $[M_b,T]$, with $T$ a Calderón-Zygmund operator, are bounded on $L^p(\Sigma,\mu)$ whenever $b\in\BMO(\Sigma,\mu)$ (cf. \cite{CoiRocWei76}). A weighted version of this in the setting of spaces of homogeneous type can be found in \cite[Theorem 2.16]{HofMitTay09}. Specifically, they prove that for each $p\in(1,\infty)$
there is a non-decreasing function $\Phi_p:[1,\infty)\to[1,\infty)$ such that
\begin{equation}\label{eq:commutator-bounded}
\| [M_b,T] \|_{L^p(\Sigma,w)\to L^p(\Sigma,w)} \leq \Phi_p([w]_{A_p}) \| b \|_{\BMO(\Sigma,\mu)}.
\end{equation}
Finally, Uchiyama proved that the commutator $[M_b,T]$ is compact in $L^p(\mathbb{R}^n,\mathcal{L}^n)$ if and only if $b \in \VMO(\mathbb{R}^n,\mathcal{L}^n)$, where $\mathcal{L}^n$ denotes the $n$-dimensional Lebesgue measure in $\mathbb{R}^n$ (cf. \cite{Uch78}).
The following theorem is a quantitative version of Uchiyama's theorem, whose proof for spaces of homogeneous type can be found in \cite[Theorem 2.19]{HofMitTay09}.
In the next result, and elsewhere in this section, $\dist_{\BMO}$ denotes the distance in the 
$\BMO(\Sigma, \mu)$ norm.

\begin{theorem} \label{thm:Uchiyama-quantitative}
Suppose that $(\Sigma, d, \mu)$ is a compact space of homogeneous type, 
and assume that $T$ is a Calderón–Zygmund operator.
Then for every $p \in (1,\infty)$, there exists $C\in(0,\infty)$ such that
\begin{equation}
\| [M_{b},T] \|^{\ess}_{L^{p}(\Sigma,\mu)\to L^{p}(\Sigma,\mu)}
\le C\, \mathrm{dist}_{\BMO}\big(b, \mathrm{VMO}(\Sigma, \mu)\big)
\end{equation}
for every $b \in \mathrm{BMO}(\Sigma, \mu)$.
In particular, if $b \in \mathrm{VMO}(\Sigma, \mu)$,
then the commutator $[M_{b}, T]$ is compact in $L^{p}(\Sigma, \mu)$ for each $p \in (1,\infty)$.
\end{theorem}

Our interpolation theorems allow us to extend this theorem to other function spaces.

\begin{theorem}
Suppose that $(\Sigma, d, \mu)$ is a compact, non-atomic space of homogeneous type, 
and assume that $T$ is a Calderón–Zygmund operator. Then, the following holds:
\begin{enumerate}
\item Let $p \in (1,\infty)$ and $w \in A_p(\Sigma)$. Then, there exist $\theta\in(0,1)$ and $C\in(0,\infty)$ such that
\begin{equation}
\| [M_{b},T] \|^{\ess}_{L^{p}(\Sigma,w)\to L^{p}(\Sigma,w)}
\leq C\, \| b \|_{\BMO(\Sigma,\mu)}^{\theta} \dist_{\BMO}\!\big(b,\VMO(\Sigma,\mu)\big)^{1-\theta}
\end{equation}
for every  $b \in \mathrm{BMO}(\Sigma, \mu)$.
In particular, if $b \in \mathrm{VMO}(\Sigma, \mu)$,
then $[M_{b}, T]$ is compact in $L^{p}(\Sigma, w)$.
\item Let $\mathbb{X}$ be a RIBFS over $(\Sigma,\mu)$ such that $1<p_{\mathbb{X}}\leq q_{\mathbb{X}} < \infty$ and fix 
$w \in A_{p_{\mathbb{X}}}(\Sigma)$ with $w(\Sigma)=\mu(\Sigma)$. Then, there exist $\theta\in(0,1)$ and $C\in(0,\infty)$ such that
\begin{equation}
\| [M_{b},T] \|^{\ess}_{\mathbb{X}(w)\to \mathbb{X}(w)}
\le C\, \| b \|_{\BMO(\Sigma,\mu)}^{\theta}  \dist_{\BMO}\!\big(b, \VMO(\Sigma, \mu)\big)^{1-\theta}
\end{equation}
for every  $b \in \mathrm{BMO}(\Sigma, \mu)$.
In particular, if $b \in \mathrm{VMO}(\Sigma, \mu)$,
then $[M_{b}, T]$ is compact in $\mathbb{X}(w)$.
\item Assume that $(\Sigma,\mu)$ is complete and let $p(\cdot)\in{\rm LH}(\Sigma)$ satisfy $1<p_{-} \leq p_+ < \infty$. Then, there exist $\theta\in(0,1)$ and $C\in(0,\infty)$ such that
\begin{equation}
\| [M_{b},T] \|^{\ess}_{L^{p(\cdot)}(\Sigma,\mu)\to L^{p(\cdot)}(\Sigma,\mu)}
\le C\, \| b \|_{\BMO(\Sigma,\mu)}^{\theta} \dist_{\BMO}\!\big(b, \mathrm{VMO}(\Sigma, \mu)\big)^{1-\theta}
\end{equation}
for every  $b \in \mathrm{BMO}(\Sigma, \mu)$.
In particular, if $b \in \mathrm{VMO}(\Sigma, \mu)$,
then $[M_{b}, T]$ is compact in $L^{p(\cdot)}(\Sigma, \mu)$.
\end{enumerate}
\end{theorem}

\begin{proof}
We begin by noting that, under the hypotheses above, whenever $b\in\BMO(\Sigma,\mu)$, the commutator $[M_b,T]$ is bounded on $L^p(\Sigma,w)$, $\mathbb{X}(w)$, and $L^{p(\cdot)}(\Sigma,\mu)$, with operator norm bounded by $C \| b \|_{\BMO(\Sigma,\mu)}$. This is a consequence of \eqref{eq:commutator-bounded} and extrapolation (cf. Theorem \ref{thm:extrapolation-BFS}, along with Remarks \ref{rem:RdF}, \ref{rem:extrapolation-X}, and \ref{rem:extrapolation-Lpvar}).

To prove $(a)$, fix $p\in(1,\infty)$ and $w\in A_p(\Sigma)$. 
As a consequence of Corollary \ref{cor:compact-interpol-1}, Theorem \ref{thm:Uchiyama-quantitative}, and \eqref{eq:commutator-bounded}, there are $\theta\in(0,1)$ and $C\in(0,\infty)$ such that $w^{1/\theta} \in A_p(\Sigma)$ and
\begin{align}
(\|[M_{b},T]\|^{\ess}_{L^{p}(\Sigma,w)\to L^{p}(\Sigma,w)}) & \leq C
(\|[M_{b},T]\|^{\ess}_{L^{p}(\Sigma,\mu)\to L^{p}(\Sigma,\mu)})^{1-\theta}
\nonumber \\ & \quad\,\,\, \cdot 
(\|[M_{b},T]\|^{\ess}_{L^{p}(\Sigma,w^{1/\theta})\to L^{p}(\Sigma,w^{1/\theta})})^{\theta}
\nonumber \\ & \leq 
C\dist_{\rm BMO}(b,\VMO(\Sigma,\mu))^{\,1-\theta}
\| b \|_{\BMO(\Sigma,\mu)}^{\theta}.
\end{align}

To address the proof of $(b)$, fix $w \in A_{p_{\mathbb{X}}}(\Sigma)$ with $w(\Sigma)=\mu(\Sigma)$. As a consequence of \eqref{eq:Ap-p2} and \eqref{eq:openness}, $w \in A_{p_{\mathbb{X}}-\varepsilon}(\Sigma)$ 
and $w \in A_{q_{\mathbb{X}}+\varepsilon}(\Sigma)$ for some $\varepsilon>0$ with 
$p_{\mathbb{X}}-\varepsilon \in(1,\infty)$. Applying Theorem \ref{thm:compact-interpol-Xw}, $(a)$, and \eqref{eq:commutator-bounded}, there are $C\in(0,\infty)$ and $\theta_0,\theta'_0,\theta_1\in(0,1)$ such that
\begin{align}
\|[M_{b},T]\|^{\ess}_{\mathbb{X}(w)\to \mathbb{X}(w)}
& \leq C \big( (\|[M_{b},T]\|^{\ess}_{L^{p_\mathbb{X}-\varepsilon}(\Sigma,w)\to L^{p_\mathbb{X}-\varepsilon}(\Sigma,w)})^{1-\theta_0} 
\nonumber \\ & \quad \,\, \cdot
(\|[M_{b},T]\|^{\ess}_{L^{q_\mathbb{X}+\varepsilon}(\Sigma,w)\to L^{q_\mathbb{X}+\varepsilon}(\Sigma,w)})^{\theta_0}
\nonumber \\ & \quad \,\, +
(\|[M_{b},T]\|^{\ess}_{L^{p_\mathbb{X}-\varepsilon}(\Sigma,w)\to L^{p_\mathbb{X}-\varepsilon}(\Sigma,w)})^{1-\theta'_0} 
\nonumber \\ & \quad \,\, \cdot
(\|[M_{b},T]\|^{\ess}_{L^{q_\mathbb{X}+\varepsilon}(\Sigma,w)\to L^{q_\mathbb{X}+\varepsilon}(\Sigma,w)})^{\theta'_0} \big)
\nonumber \\ & \leq C  \Big(
\big(  \| b \|_{\BMO(\Sigma,\mu)}^{\theta_1} \dist_{\rm BMO}(b,\VMO(\Sigma,\mu))^{\,1-\theta_1} \big)^{1-\theta_0}
\| b \|_{\BMO(\Sigma,\mu)}^{\theta_0}
\nonumber \\ & \quad \,\, +
\big(  \| b \|_{\BMO(\Sigma,\mu)}^{\theta_1} \dist_{\rm BMO}(b,\VMO(\Sigma,\mu))^{\,1-\theta_1} \big)^{1-\theta'_0}
\| b \|_{\BMO(\Sigma,\mu)}^{\theta'_0}
\nonumber \\ & = C  \big(
\| b \|_{\BMO(\Sigma,\mu)}^{\theta_1(1-\theta_0)+\theta_0} \dist_{\rm BMO}(b,\VMO(\Sigma,\mu))^{(1-\theta_1)(1-\theta_0)}
\nonumber \\ & \quad \,\, + 
\| b \|_{\BMO(\Sigma,\mu)}^{\theta_1(1-\theta'_0)+\theta'_0} \dist_{\rm BMO}(b,\VMO(\Sigma,\mu))^{(1-\theta_1)(1-\theta'_0)} \big)
\nonumber \\ & \leq C
\| b \|_{\BMO(\Sigma,\mu)}^{\theta} \dist_{\rm BMO}(b,\VMO(\Sigma,\mu))^{1-\theta},
\end{align}
where $\theta:=\max\{ \theta_1(1-\theta_0)+\theta_0,\,\theta_1(1-\theta'_0)+\theta'_0 \}$ and we have used that for every $b\in \BMO(\Sigma,\mu)$ one has $\dist_{\rm BMO}(b,\VMO(\Sigma,\mu)) \leq \| b \|_{\BMO(\Sigma,\mu)}$.

Finally, to prove the result for variable exponent Lebesgue spaces described in $(c)$, fix $p(\cdot) \in {\rm LH}(\Sigma)$ 
satisfying $1<p_{-}\leq p_{+}<\infty$, $p_1\in(1,\infty)$, and apply Corollary \ref{cor:compact-interpol-Lpvar-2} to deduce that 
there are $q(\cdot) \in {\rm LH}(\Sigma)$ with  $1<q_{-}\leq q_{+} < \infty$, $\theta\in(0,1)$ and $C\in(0,\infty)$ such that
\begin{align}
\| [M_{b},T] \|^{\ess}_{L^{p(\cdot)}(\Sigma,\mu)\to L^{p(\cdot)}(\Sigma,\mu)}
& \leq  C (\| [M_{b},T] \|^{\ess}_{L^{p_1}(\Sigma,\mu)\to L^{p_1}(\Sigma,\mu)})^{1-\theta}
\nonumber\\ & \qquad 
\cdot (\| [M_{b},T] \|^{\ess}_{L^{q(\cdot)}(\Sigma,\mu)\to L^{q(\cdot)}(\Sigma,\mu)})^{\theta} 
\nonumber\\ & \leq 
C\dist_{\rm BMO}(b,\VMO(\Sigma,\mu))^{1-\theta}
\| b \|_{\BMO(\Sigma,\mu)}^{\theta},
\end{align}
where the second inequality follows from Theorem \ref{thm:Uchiyama-quantitative}, \eqref{eq:commutator-bounded}, and extrapolation.
\end{proof}

We close this section by noting that other results of this type are proved in \cite[Section 2.7]{MitMitMit-book3} using different techniques whenever $\Sigma$ is a compact subset of $\mathbb{R}^n$ satisfying certain geometric assumptions.

\subsection{Fredholmness and invertibility of operators} \label{sec:Fred}

Given two Banach spaces $X$ and $Y$, recall that a bounded linear operator $T:X \to Y$
is called a {\tt Fredholm} {\tt operator} if $\Im(T)$ is closed and $\dim(\Ker(T))$ and $\dim(\Coker(T))$ 
are both finite, where $\Ker(T)$ stands for the kernel of $T$, $\Coker(T):= Y/\Im(T)$ is the 
cokernel and $\Im(T)$ stands for the range. In this case, the {\tt Fredholm} {\tt index}
of the operator is defined as
\begin{equation}
\ind(T):=\dim(\Ker(T)) - \dim(\Coker(T)).
\end{equation}

We first gather some well-known properties of Fredholm operators in the following lemma, 
whose proof can be found in \cite[Chapter 5]{Sch02} or \cite{Con99}.

\begin{lemma} \label{lem:Fredholm}
Let $X,Y$ be Banach spaces and let $T:X \to Y$ be a Fredholm operator. Then,
\begin{enumerate}
\item \label{lem:Fredholm-a} If $K:X \to Y$ is compact, then $T+K$ is a Fredholm operator and $\ind(T+K)=\ind(T)$.
\item \label{lem:Fredholm-b} If $Z$ is a Banach space and $S:Y \to Z$ is another Fredholm operator,
then $S \circ T$ is also a Fredholm operator and $\ind(S \circ T)=\ind(S) + \ind(T)$.
\item \label{lem:Fredholm-d} There exists $\varepsilon_{T}>0$ depending on $T$ such that for any bounded linear operator $S:X \to Y$
with $\norm{S}_{X \to Y}<\varepsilon_{T}$, one has that $T+S$ is a Fredholm operator and $\ind(T+S)=\ind(T)$.
\end{enumerate}
\end{lemma}

The next proposition extends Lemma \ref{lem:Fredholm}\ref{lem:Fredholm-a}
replacing the compact operator $K$ with an operator $S$ which is close to the
subspace of compact operators.

\begin{proposition}\label{prop:distFred}
Let $X,Y$ be Banach spaces and $T:X \to Y$ a Fredholm operator.
Then, there exists $\varepsilon_T>0$ such that
for every bounded linear operator $S:X \to Y$ with $\| S \|^{\ess}_{X \to Y}<\varepsilon_T$,
\begin{equation}\label{eq:invertiblecond-1}
\parbox{10cm}{\begin{center}
$T+S$ is a Fredholm operator and $\ind(T+S)=\ind(T)$.
\end{center}}
\end{equation}
If we additionally assume that $\ind(T)=0$ and $T+S:X\to Y$ is injective, then $T+S$ is invertible.

\end{proposition}
\begin{proof}
Assume that $\| S \|^{\ess}_{X \to Y}<\varepsilon_T$, with $\varepsilon_T>0$ as in Lemma \ref{lem:Fredholm}\ref{lem:Fredholm-d}. Then, there exists a compact operator $K:X \to Y$ such that
$\norm{S-K}_{X\to Y}<\varepsilon_T$.
We now apply Lemma \ref{lem:Fredholm}\ref{lem:Fredholm-a} 
and Lemma \ref{lem:Fredholm}\ref{lem:Fredholm-d}, which guarantee that
$T+S=T+(S-K)+K$ is a Fredholm operator and $\ind(T+S)=\ind(T)$.
This proves \eqref{eq:invertiblecond-1}. 
The proof of the last claim follows immediately, since
any injective Fredholm operator with index $0$ is invertible.
\end{proof}

In conjunction with the results presented in Section \ref{sec:interpol}, this can be used to obtain Fredholmness and invertibility for operators acting on weighted Lebesgue spaces, weighted rearrangement invariant Banach function spaces, and variable exponent Lebesgue spaces.
One relevant example is the following result, where we establish Fredholmness and invertibility for operators of the form $\lambda I + S$, where $I$ is the identity operator, thereby capturing the spectral properties of $S$.

\begin{proposition}\label{prop:Fred-S}
Let $(\Sigma,d,\mu)$ be a non-atomic space of homogeneous type and fix two exponents $p_0,p_1\in(1,\infty)$. Assume that $S$ is a linear operator such that for every weight $w_0 \in A_{p_0}(\Sigma)$ one has
\begin{equation}\label{eq:Fred-hyp}
\|Sf\|_{L^{p_0}(\Sigma,w_0)} \leq \Phi([w_0]_{A_{p_0}})\|f\|_{L^{p_0}(\Sigma,w_0)}
\quad \forall f \in L^{p_0}(\Sigma,w_0), 
\end{equation}
where $\Phi: [1, \infty) \to [1, \infty)$ is a non-decreasing function.
Then, the following holds:
\begin{enumerate}
\item For every $\lambda_0>0$, $p\in(1,\infty)$, and $w\in A_p(\Sigma)$ there exists $\varepsilon>0$ such that if $\| S \|^{\ess}_{L^{p_1}(\Sigma,\mu) \to L^{p_1}(\Sigma,\mu)}<\varepsilon$,
\begin{equation}
\parbox{8cm}{\begin{center}\label{eq:Fred-sp1}
$\lambda I+S$ is a Fredholm operator with index $0$ on $L^{p}(\Sigma,w)$ whenever $|\lambda|>\lambda_0$.
\end{center}}
\end{equation}
\item Let $\mathbb{X}$ be a RIBFS over $(\Sigma,\mu)$ with $1<p_{\mathbb{X}}\leq q_{\mathbb{X}} < \infty$ and fix $w \in A_{p_{\mathbb{X}}}(\Sigma)$ with $w(\Sigma)=\mu(\Sigma)$.
Then, for every $\lambda_0>0$, there exists $\varepsilon>0$ such that if $\| S \|^{\ess}_{L^{p_1}(\Sigma,\mu) \to L^{p_1}(\Sigma,\mu)}<\varepsilon$,
\begin{equation}\label{eq:Fred-sp2}
\parbox{8cm}{\begin{center}
$\lambda I+S$ is a Fredholm operator with index $0$ on $\mathbb{X}(w)$ whenever $|\lambda|>\lambda_0$.
\end{center}}
\end{equation}
\item Assume that $(\Sigma,\mu)$ is complete and
let $p(\cdot)\in{\rm LH}(\Sigma)$ satisfy $1<p_{-} \leq p_+ < \infty$.
Then, for every $\lambda_0>0$, there exists $\varepsilon>0$ such that if $\| S \|^{\ess}_{L^{p_1}(\Sigma,\mu) \to L^{p_1}(\Sigma,\mu)}<\varepsilon$,
\begin{equation}\label{eq:Fred-sp3}
\parbox{8cm}{\begin{center}
$\lambda I+S$ is a Fredholm operator with index $0$ on $L^{p(\cdot)}(\Sigma,\mu)$ whenever $|\lambda|>\lambda_0$.
\end{center}}
\end{equation}
\item \label{prop:Fred-S-d} Under the additional assumptions that $(\Sigma,d,\mu)$ is bounded and there is $\varepsilon>0$ such that $\lambda I+S$ is injective on $L^q(\Sigma,\mu)$ for every $q\in(1,1+\varepsilon)$, we have that
the operators $\lambda I+S$ in \eqref{eq:Fred-sp1}-\eqref{eq:Fred-sp3} are actually invertible.
\end{enumerate}
\end{proposition}
\begin{proof}
We start by noting that, under the hypotheses above, $S$ (and hence $\lambda I + S$ for all $\lambda\in\mathbb{R}$) is bounded on $L^p(\Sigma,w)$, $\mathbb{X}(w)$, and $L^{p(\cdot)}(\Sigma,\mu)$ as a consequence of \eqref{eq:Fred-hyp} and extrapolation (cf. Theorem \ref{thm:extrapolation-BFS}, along with Remarks \ref{rem:RdF}, \ref{rem:extrapolation-X}, and \ref{rem:extrapolation-Lpvar}).

To address the proof of $(a)$, fix $p\in(1,\infty)$, $w\in A_p(\Sigma)$, and $\lambda_0>0$. As a consequence of Proposition~\ref{prop:distFred}, there is $\varepsilon_0>0$ such that if
$\|\frac{1}{\lambda}S \|^{\ess}_{L^p(\Sigma,w)\to L^p(\Sigma,w)} < \varepsilon_0$, then
\begin{equation}\label{eq:Fred.inproof-2}
\parbox{14cm}{\begin{center}
$\displaystyle I+\frac{1}{\lambda}S$ is a Fredholm operator on $L^p(\Sigma,w)$ and $\displaystyle \ind\Big(I+\frac{1}{\lambda}S\Big)=0$.
\end{center}}
\end{equation}
On the other hand, Corollary \ref{cor:compact-interpol-2} guarantees that there are $\theta\in(0,1)$ and  
$C_0\in[1,\infty)$ such that
\begin{equation}\label{eq:Fred.inproof-1}
\|S\|^{\ess}_{L^p(\Sigma,w)\to L^p(\Sigma,w)} \leq C_0
(\|S\|^{\ess}_{L^{p_1}(\Sigma,\mu)\to L^{p_1}(\Sigma,\mu)})^{\theta}.
\end{equation}
We now take $\varepsilon < (\lambda_0 \varepsilon_0 / C_0 )^{1/\theta}$, where 
$C_0$ is the constant from \eqref{eq:Fred.inproof-1},
and assume that $\| S \|^{\ess}_{L^{p_1}(\Sigma,\mu) \to L^{p_1}(\Sigma,\mu)}<\varepsilon$.
Then, whenever $|\lambda|>\lambda_0$, \eqref{eq:Fred.inproof-1} yields
\begin{equation}
\Big\|\frac{1}{\lambda}S \Big\|^{\ess}_{L^p(\Sigma,w)\to L^p(\Sigma,w)} 
\leq \frac{C_0}{|\lambda|}  \varepsilon^{\theta}
< \varepsilon_0.
\end{equation}
As a consequence of \eqref{eq:Fred.inproof-2},
$I + \frac{1}{\lambda}S$, and hence
$\lambda I + S = \lambda\big( I + \frac{1}{\lambda}S \big)$, is a Fredholm
operator with index $0$ on $L^p(\Sigma,w)$.

We now proceed to prove the claim in $(b)$. Fix $w\in A_{p_\mathbb{X}}(\Sigma)$ with $w(\Sigma)=\mu(\Sigma)$ and $\lambda_0>0$.
As a consequence of Proposition~\ref{prop:distFred}, there is a threshold $\varepsilon_0>0$ such that if
$\|\frac{1}{\lambda}S \|^{\ess}_{\mathbb{X}(w)\to \mathbb{X}(w)} < \varepsilon_0$, then
\begin{equation}\label{eq:Fred.inproof-4}
\parbox{14cm}{\begin{center}
$\displaystyle I+\frac{1}{\lambda}S$ is a Fredholm operator on $\mathbb{X}(w)$ and $\displaystyle \ind\Big(I+\frac{1}{\lambda}S\Big)=0$.
\end{center}}
\end{equation}
On the other hand, using Corollary \ref{cor:compact-interpol-Xw-2}, there exist $\theta,\theta'\in(0,1)$ and $C_0\in(0,\infty)$ such that
\begin{equation}\label{eq:Fred.inproof-3}
\|S\|^{\ess}_{\mathbb{X}(w)\to \mathbb{X}(w)}
\leq C_0 \big(
(\|S\|^{\ess}_{L^{p_1}(\Sigma,\mu)\to L^{p_1}(\Sigma,\mu)})^{\theta}
+
(\|S\|^{\ess}_{L^{p_1}(\Sigma,\mu)\to L^{p_1}(\Sigma,\mu)})^{\theta'} \big).
\end{equation}
Take $\varepsilon>0$ small enough so that 
$C_0 (\varepsilon^\theta + \varepsilon^{\theta'} ) / \lambda_0 < \varepsilon_0$,
where $C_0$ is the constant from \eqref{eq:Fred.inproof-3},
and assume that $\| S \|^{\ess}_{L^{p_1}(\Sigma,\mu) \to L^{p_1}(\Sigma,\mu)}<\varepsilon$.
Then, whenever $|\lambda|>\lambda_0$, \eqref{eq:Fred.inproof-3} gives
\begin{equation}
\Big\|\frac{1}{\lambda}S \Big\|^{\ess}_{\mathbb{X}(w)\to \mathbb{X}(w)} 
\leq \frac{C_0}{|\lambda|}  \big(\varepsilon^{\theta} + \varepsilon^{\theta'}\big)
< \varepsilon_0.
\end{equation}
As a consequence of \eqref{eq:Fred.inproof-4},
$I + \frac{1}{\lambda}S$, and hence
$\lambda I + S = \lambda\big( I + \frac{1}{\lambda}S \big)$, is a Fredholm
operator with index $0$ on $\mathbb{X}(w)$.

To prove the result in $(c)$, fix $p(\cdot) \in {\rm LH}(\Sigma)$ with $1<p_{-}\leq p_{+}<\infty$ and $\lambda_0>0$.
As a consequence of Corollary \ref{cor:compact-interpol-Lpvar-2},
there are $q(\cdot) \in {\rm LH}(\Sigma)$ verifying $1<q_{-}\leq q_{+} < \infty$, $\theta\in(0,1)$ and $C,C_0\in(0,\infty)$ such that
\begin{align}
\| S \|^{\ess}_{L^{p(\cdot)}(\Sigma,\mu)\to L^{p(\cdot)}(\Sigma,\mu)}
& \leq   C (\| S \|^{\ess}_{L^{p_1}(\Sigma,\mu)\to L^{p_1}(\Sigma,\mu)})^{1-\theta}
(\| S \|^{\ess}_{L^{q(\cdot)}(\Sigma,\mu)\to L^{q(\cdot)}(\Sigma,\mu)})^{\theta}
\nonumber \\ & \leq  C_0 (\| S \|^{\ess}_{L^{p_1}(\Sigma,\mu)\to L^{p_1}(\Sigma,\mu)})^{1-\theta},
\end{align}
where the second inequality follows from \eqref{eq:Fred-hyp} and extrapolation.
The remainder of the proof of $(c)$ proceeds analogously to the proofs of $(a)$ and $(b)$.

Finally, to prove $(d)$, we start by fixing $p\in(1,\infty)$ and $w\in A_p(\Sigma)$. As a consequence of \eqref{eq:Ap-p2}-\eqref{eq:openness}, there exists $\delta>0$ such that $w\in A_{p/q}(\Sigma,\mu)$ for every $q\in(1,1+\delta)$. Then, for every
$f\in L^p(\Sigma,w)$, applying H\"older's inequality, we obtain 
\begin{align}\label{eq:inclusion-1}
\int_{\Sigma}\left|f\right|^q\,d\mu
&=\int_{\Sigma}\left|f\right|^q w^{q/p}w^{-q/p}\,d\mu 
\leq\left(\int_{\Sigma}\left| f\right|^{p}dw\right)^{q/p}
\left(\int_{\Sigma}w^{(-q/p)\left( p/q\right)'}d\mu\right)^{1/(p/q)'}
\nonumber\\[4pt]
&\leq\norm{f}_{L^p(\Sigma,w)}^{q}\left(\int_{\Sigma}w^{1-(p/q)'}d\mu\right)^{1/(p/q)'}<\infty,
\end{align}
where the last term is finite because $w\in A_{p/q}(\Sigma,\mu)$ implies that
$w^{1-(p/q)'}\in L^1(\Sigma)$ since we are assuming that $(\Sigma,d,\mu)$ is bounded (see \eqref{eq:Ap-p1}).
This proves $L^p(\Sigma,w)\subseteq L^{q}(\Sigma,\mu)$ 
whenever $q\in(1,1+\delta)$. Since we are assuming that $\lambda I + S$ is injective on $L^{q}(\Sigma,\mu)$ for some $q\in(1,1+\delta)$, we deduce that $\lambda I + S$ is injective on $L^p(\Sigma,w)$ as well and the result follows at once, since any injective Fredholm operator with index $0$ is invertible.
For the case of RIBFS, note first that $\mu(\Sigma)<+\infty$ whenever $(\Sigma,d,\mu)$ is bounded (see \eqref{eq:doubling}). Hence, from \eqref{eq:Xw-inclusion}, we have that $\mathbb{X}(w)\subseteq L^{p}(\Sigma,w)$ whenever $p\in(1,p_{\mathbb{X}})$ and the result follows as before.
The case of variable exponent Lebesgue spaces is handled similarly since
$L^{p(\cdot)}(\Sigma,\mu) \subseteq L^{p_{-}}(\Sigma,\mu)$  under the current hypotheses (cf. \cite[Theorem 3.3.11]{DieHarHasRuz11}).
\end{proof}

\subsection{Essential norm of Chord-Dot-Normal singular integral operators}

The aim of this section is to bound the essential norm in certain function spaces of the so-called \textit{Chord–Dot–Normal} singular integral operators, whose defining feature is that their kernel is given by an inner product with the outward unit normal vector.
This class of operators includes harmonic double layer potentials, among others (cf. Remark \ref{rem:RK} below).
While this is covered in \cite[Chapter 5]{MitMitMit-book5}, we provide different proofs via our compact interpolation approach.

We start by introducing several notions from geometric measure theory. 
Recall that $\mathcal{L}^{n}$ denotes the $n$-dimensional Lebesgue measure in $\mathbb{R}^n$.
In what follows, we fix $n\geq 2$, we denote the $(n-1)$-dimensional Hausdorff measure in $\mathbb{R}^n$ as $\mathcal{H}^{n-1}$, and we let $\Omega$ be an $\mathcal{L}^n$-measurable, open, nonempty proper subset of $\mathbb{R}^n$.
We define the {\tt measure theoretic boundary} of $\Omega$ as
\begin{equation}
\partial_{*}\Omega := \Big\{x\in\partial\Omega:\,\limsup_{r\to 0^{+}}
\frac{\mathcal{L}^{n}(B(x,r)\cap\Omega)}{r^{n}}> 0, \,\,\,
\limsup_{r\to 0^{+}}\frac{\mathcal{L}^{n}(B(x,r)\setminus\Omega)}{r^{n}}>0\Big\},
\end{equation}
and we say that $\Omega$ has {\tt locally finite perimeter} if 
\begin{equation}
{\mathcal{H}}^{n-1}\big(\partial_{*}\Omega\cap K) <+\infty \quad
\text{ for every compact set }\,K\subset{\mathbb{R}}^n.
\end{equation}
If $\Omega$ has locally finite perimeter, one can deduce (cf., e.g.,  \cite[Chapter 5]{EvaGar92} 
for details) that there exists a $\mathbb{R}^{n}$-valued measurable 
function $\nu$ on $\partial_{*}\Omega$ such that $|\nu(x)|=1$ $\mathcal{H}^{n-1}$-a.e. on $\partial_{*}\Omega$ and 
\begin{equation}
\int_{\Omega}\div\vec{F}\,d\mathcal{L}^{n}
= \int_{\partial_\ast\Omega}\nu\cdot\big(\vec{F}\big|_{\partial\Omega}\big)\,d{\mathcal{H}}^{n-1}
\quad\text{ for every }\vec{F}\in\big[\mathscr{C}_0^1(\mathbb{R}^{n})\big]^n.
\end{equation}
The function $\nu$ is called {\tt geometric} {\tt measure} {\tt theoretic} {\tt outward} 
{\tt unit} {\tt normal} to $\Omega$.

\begin{definition}\label{def:AR}
A closed set $\Sigma\subset{\mathbb{R}}^{n}$ is called {\tt Ahlfors} {\tt regular} {\tt set} if there is a constant $C\in[1,\infty)$ such that 
\begin{equation}\label{eq:AR}
r^{n-1}/C\leq{\mathcal{H}}^{n-1}\big(B(x,r)\cap\Sigma\big)\leq Cr^{n-1},
\quad\forall\,r\in\big(0,2\,{\rm diam}(\Sigma)\big),\quad\forall\,x\in\Sigma. 
\end{equation}
\end{definition}

\begin{definition}
An open, nonempty, proper subset $\Omega$ of $\mathbb{R}^n$ is called an {\tt Ahlfors regular domain} if $\partial\Omega$ is an Ahlfors regular set and
$\mathcal{H}^{n-1}(\partial\Omega\setminus\partial_{*}\Omega)=0$.
\end{definition}

If $\Omega\subset\mathbb{R}^n$ is an Ahlfors regular domain, then $\sigma:={\mathcal{H}}^{n-1}\lfloor\partial\Omega$ is a complete, Radon, doubling measure on $\partial\Omega$ (cf. \cite[(5.9.14)]{MitMitMit-book1}), which we will refer to as the {\tt surface measure} of $\Omega$.

\begin{definition}\label{def:john}
An open, nonempty, proper subset $\Omega$ of ${\mathbb{R}}^{n}$ is said to satisfy a 
{\tt local} {\tt John} {\tt condition} if there are $\theta\in(0,1)$ and $R>0$
{\rm (}with the requirement that $R=\infty$ whenever $\partial\Omega$ is unbounded{\rm )} 
such that for every $x\in\partial\Omega$ and $r\in(0,R)$ there exists a point $x_r\in B(x,r)\cap\Omega$ 
such that $B(x_r,\theta r)\subseteq\Omega$ and with the property that for each
$y\in B(x,r)\cap\partial\Omega$ there is a rectifiable path
$\gamma_{y}:[0,1]\to\overline{\Omega}$ whose length is $\leq\theta^{-1}r$ and such that 
\begin{equation}\label{eq:john}
\gamma_{y}(0)=y,\quad\gamma_{y}(1)=x_r,\quad
{\rm dist}(\gamma_{y}(t)\,,\,\partial\Omega)>\theta|\gamma_{y}(t)-y|\,\text{ for all }
\,t\in(0,1]. 
\end{equation}
Finally, an open, nonempty set $\Omega\subset{\mathbb{R}}^{n}$ which is not dense in
${\mathbb{R}}^n$ is said to satisfy a {\tt two}-{\tt sided} {\tt local} {\tt John}
{\tt condition} if both $\Omega$ and ${\mathbb{R}}^n\setminus\overline{\Omega}$ satisfy a local 
John condition.
\end{definition}

The following theorem is one of the main results in \cite[Theorem 4.36]{HofMitTay09} and it bounds the essential norm of Chord-Dot-Normal operators on Lebesgue spaces.
As usual, $\dist_{\BMO}$ denotes the distance in the 
$\BMO(\partial\Omega, \sigma)$ norm.

\begin{theorem}\label{thm:SIOCompact-unweighted}
Let $\Omega\subset\mathbb{R}^{n}$ be an Ahlfors regular domain satisfying a two-sided local 
John condition with compact boundary. Abbreviate 
$\sigma:={\mathcal{H}}^{n-1}\lfloor\partial\Omega$, fix 
$p\in(1,\infty)$, and consider a sufficiently large integer $N=N(n)\in\mathbb{N}$.
Given a function
$k\in\mathscr{C}^{N}(\mathbb{R}^{n}\setminus\{0\})$ such that 
\begin{equation} \label{eq:SIOCompact-khyp}
k(x)=k(-x), \quad k(\lambda x)=\lambda^{-n}k(x) \quad \text{ for every }\lambda>0 \text{ and } x\in\mathbb{R}^{n}\setminus\{0\},
\end{equation}
consider the operator
\begin{equation}\label{eq:SIOCompact-unweighted}
Tf(x)=\lim_{\eta\to0^{+}}
\int_{\substack{y\in\partial\Omega\\|x-y|>\eta}}\langle x-y,\nu(y)\rangle k(x-y)f(y)d\sigma(y),
\quad x\in\partial\Omega.
\end{equation}
Then, $T$ is a well-defined linear bounded operator on $L^p(\partial\Omega,\sigma)$.
Moreover, for every $\varepsilon>0$ there exists $\delta>0$
such that
\begin{equation}\label{eq:SIOCompact-1}
\dist_{\BMO}\!\big(\nu,[\VMO(\partial\Omega,\sigma)]^n\big)<\delta
\quad \text{ implies } \quad
\|T\|^{\ess}_{L^p(\partial\Omega,\sigma)\to L^p(\partial\Omega,\sigma)}<\varepsilon.
\end{equation}
In particular, if $\nu\in[\VMO(\partial\Omega,\sigma)]^n$ then $T:L^p(\partial\Omega,\sigma)\longrightarrow L^p(\partial\Omega,\sigma)$ is compact.
\end{theorem}

The operators $T$ defined in \eqref{eq:SIOCompact-khyp}-\eqref{eq:SIOCompact-unweighted} are actually bounded on $L^p(\partial\Omega,w)$ for all $p\in(1,\infty)$ and  $w\in A_p(\partial\Omega)$ (cf. \cite[Proposition 3.4]{MarMarMitMitMit22}). Specifically, under the hypotheses of Theorem \ref{thm:SIOCompact-unweighted}, there is a non-decreasing function $\Phi_p: [1, \infty) \to [1, \infty)$ (depending only on $n$, $p$, $\norm{k}_{\mathscr{C}^N(S^{n-1})}$, and the geometric constants of $\partial\Omega$) such that
\begin{equation}\label{eq:T-bounded-Lpw}
\|Tf\|_{L^p(\partial\Omega,w)} \leq \Phi_p([w]_{A_p}) \|f\|_{L^p(\partial\Omega,w)}
\qquad \forall f \in L^p(\partial\Omega,w).
\end{equation}
This leads us to the next theorem, in which we extend Theorem \ref{thm:SIOCompact-unweighted} to additional function spaces.

\begin{theorem}\label{thm:SIOCompact}
Let $\Omega\subset\mathbb{R}^{n}$ be an Ahlfors regular domain satisfying a two-sided local 
John condition with compact boundary. Abbreviate 
$\sigma:={\mathcal{H}}^{n-1}\lfloor\partial\Omega$ and consider a sufficiently large integer $N=N(n)\in\mathbb{N}$.
Given a function
$k\in\mathscr{C}^{N}(\mathbb{R}^{n}\setminus\{0\})$ such that 
\begin{equation}
k(x)=k(-x), \quad k(\lambda x)=\lambda^{-n}k(x) \quad \text{ for every }\lambda>0 \text{ and } x\in\mathbb{R}^{n}\setminus\{0\},
\end{equation}
consider the operator
\begin{equation}\label{thm:SIOCompact2}
Tf(x)=\lim_{\eta\to0^{+}}
\int_{\substack{y\in\partial\Omega\\|x-y|>\eta}}\langle x-y,\nu(y)\rangle k(x-y)f(y)d\sigma(y),
\quad x\in\partial\Omega.
\end{equation}
Then, the following holds:
\begin{enumerate}
\item Let $p\in(1,\infty)$ and $w \in A_p(\partial\Omega)$. Then, for every $\varepsilon>0$ there exists $\delta>0$
such that
\begin{equation}\label{eq:SIOCompact2-1}
\dist_{\BMO}\!\big(\nu,[\VMO(\partial\Omega,\sigma)]^n\big)<\delta
\quad \text{ implies } \quad
\|T\|^{\ess}_{L^p(\partial\Omega,w)\to L^p(\partial\Omega,w)}<\varepsilon.
\end{equation}
In particular, if $\nu\in[\VMO(\partial\Omega,\sigma)]^n$ then $T:L^p(\partial\Omega,w)\longrightarrow L^p(\partial\Omega,w)$ is compact.
\item Let $\mathbb{X}$ be a RIBFS over $(\partial\Omega,\sigma)$ such that $1<p_{\mathbb{X}}\leq q_{\mathbb{X}} < \infty$ and fix a weight
$w \in A_{p_{\mathbb{X}}}(\partial \Omega)$ with $w(\partial\Omega)=\sigma(\partial\Omega)$. Then, for every $\varepsilon>0$ there exists $\delta>0$
such that
\begin{equation}\label{eq:SIOCompact2-5}
\dist_{\BMO}\!\big(\nu,[\VMO(\partial\Omega,\sigma)]^n\big)<\delta
\quad \text{ implies } \quad
\|T\|^{\ess}_{\mathbb{X}(w)\to \mathbb{X}(w)}<\varepsilon.
\end{equation}
In particular, if $\nu\in[\VMO(\partial\Omega,\sigma)]^n$ then $T:\mathbb{X}(w)\longrightarrow \mathbb{X}(w)$ is compact.
\item Let $p(\cdot)\in{\rm LH}(\partial\Omega)$ satisfy $1<p_{-} \leq p_+ < \infty$. Then, for every $\varepsilon>0$ there exists $\delta>0$
such that
\begin{equation}\label{eq:SIOCompact2-3}
\dist_{\BMO}\!\big(\nu,[\VMO(\partial\Omega,\sigma)]^n\big)<\delta
\quad \text{ implies } \quad
\|T\|^{\ess}_{L^{p(\cdot)}(\partial\Omega,\sigma)\to L^{p(\cdot)}(\partial\Omega,\sigma)}<\varepsilon.
\end{equation}
In particular, if $\nu\in[\VMO(\partial\Omega,\sigma)]^n$ then $T:L^{p(\cdot)}(\partial\Omega,\sigma)\longrightarrow L^{p(\cdot)}(\partial\Omega,\sigma)$ is compact.
\end{enumerate}
\end{theorem}
\begin{proof}
First, we note that, under the stated assumptions, $T$ is bounded on $L^p(\partial\Omega,w)$, $\mathbb{X}(w)$, and $L^{p(\cdot)}(\partial\Omega,\sigma)$ as a consequence of \eqref{eq:T-bounded-Lpw} and extrapolation (cf. Theorem \ref{thm:extrapolation-BFS}, along with Remarks \ref{rem:RdF}, \ref{rem:extrapolation-X}, and \ref{rem:extrapolation-Lpvar}).

To prove $(a)$, fix $p\in(1,\infty)$, $w\in A_p(\partial\Omega)$, and $\varepsilon>0$. It follows from Corollary \ref{cor:compact-interpol-2} that there are $\theta\in(0,1)$ and
$C_0\in[1,\infty)$ such that
\begin{equation}\label{eq:SIO.inproof-1}
\|T\|^{\ess}_{L^p(\partial\Omega,w)\to L^p(\partial\Omega,w)}  \leq C_0
(\|T\|^{\ess}_{L^{p}(\partial\Omega,\sigma)\to L^{p}(\partial\Omega,\sigma)})^{\theta}.
\end{equation}
On the other hand, Theorem \ref{thm:SIOCompact-unweighted} guarantees that there exists $\delta>0$ such that if $\dist_{\BMO}\!\big(\nu,[\VMO(\partial\Omega,\sigma)]^n\big)<\delta$, then
\begin{equation}\label{eq:SIO.inproof-2}
\|T\|^{\ess}_{L^p(\partial\Omega,\sigma)\to L^p(\partial\Omega,\sigma)}<\left( \frac{\varepsilon}{C_0} \right)^{1/\theta},
\end{equation}
where $C_0$ is the constant from \eqref{eq:SIO.inproof-1}.
Combining \eqref{eq:SIO.inproof-1} and \eqref{eq:SIO.inproof-2} gives 
the result.

To prove $(b)$, fix $w\in A_{p_\mathbb{X}}(\partial\Omega)$ with $w(\partial\Omega)=\sigma(\partial\Omega)$ and $\varepsilon>0$.
Applying Corollary \ref{cor:compact-interpol-Xw-2}, for any $p_1\in(1,\infty)$ there are $\theta,\theta'\in(0,1)$ and $C_0\in(0,\infty)$ such that
\begin{equation}\label{eq:SIO.inproof-3}
\|T\|^{\ess}_{\mathbb{X}(w)\to \mathbb{X}(w)}
\leq C_0 \big( 
(\|T\|^{\ess}_{L^{p_1}(\partial\Omega,\sigma)\to L^{p_1}(\partial\Omega,\sigma)})^{\theta}
+ (\|T\|^{\ess}_{L^{p_1}(\partial\Omega,\sigma)\to L^{p_1}(\partial\Omega,\sigma)})^{\theta'} \big).
\end{equation}
The result then follows from Theorem \ref{thm:SIOCompact-unweighted} as before, since we now can take $\delta>0$ small enough so that if $\dist_{\BMO}\!\big(\nu,[\VMO(\partial\Omega,\sigma)]^n\big)<\delta$, then the right-hand side of \eqref{eq:SIO.inproof-3} is bounded by $\varepsilon$.

Finally, we address the proof of $(c)$. Fix a variable exponent $p(\cdot) \in {\rm LH}(\partial\Omega)$ 
satisfying $1<p_{-}\leq p_{+}<\infty$, $p_1\in(1,\infty)$, and $\varepsilon>0$. From Corollary \ref{cor:compact-interpol-Lpvar-2}, we deduce that 
there are $q(\cdot) \in {\rm LH}(\partial\Omega)$ with  $1<q_{-}\leq q_{+} < \infty$, $\theta\in(0,1)$ and $C,C_0\in(0,\infty)$ such that
\begin{align}\label{eq:SIO.inproof-5}
\| T \|^{\ess}_{L^{p(\cdot)}(\partial\Omega,\sigma)\to L^{p(\cdot)}(\partial\Omega,\sigma)}
& \leq  C (\| T \|^{\ess}_{L^{p_1}(\partial\Omega,\sigma)\to L^{p_1}(\partial\Omega,\sigma)})^{1-\theta}
(\| T \|^{\ess}_{L^{q(\cdot)}(\partial\Omega,\sigma)\to L^{q(\cdot)}(\partial\Omega,\sigma)})^{\theta} 
\nonumber \\ & \leq C_0
(\| T \|^{\ess}_{L^{p_1}(\partial\Omega,\sigma)\to L^{p_1}(\partial\Omega,\sigma)})^{1-\theta},
\end{align}
where the second inequality follows from \eqref{eq:T-bounded-Lpw} and extrapolation.
Once again, the conclusion follows by applying Theorem \ref{thm:SIOCompact-unweighted}.
\end{proof}

\begin{remark}\label{rem:RK}
The harmonic double layer potential is an example of an operator that falls within the scope of the previous result. 
Specifically, let $\Omega\subset\mathbb{R}^n$ be a domain as above, $\nu$ its outward unit normal vector, and $\sigma:={\mathcal{H}}^{n-1}\lfloor\partial\Omega$. Then,
the harmonic double layer potential operator is given by
\begin{equation}\label{eq:def-Kharmonic}
K_{\Delta}f(x) := \lim_{\varepsilon \to 0^+}
\frac{1}{\omega_{n-1}}
\int_{\substack{y\in\partial\Omega \\ |x - y| > \varepsilon}}
\frac{\langle \nu(y),y-x \rangle}{|x-y|^{n}} \, f(y)\, d\sigma(y),
\qquad x \in \partial\Omega,
\end{equation}
where $\omega_{n-1}$ denotes the surface measure of the unit sphere in $\mathbb{R}^n$.
The fact that the previous theorem applies to $K_\Delta$, together with the invertibility results from Section \ref{sec:Fred}, allows us to extend well-posedness of boundary value problems beyond the unweighted Lebesgue setting. This is presented in the following section.
\end{remark}

\subsection{The Dirichlet problem in non-smooth domains}
This section is devoted to providing an example of the scope of our results by establishing well-posedness for the harmonic Dirichlet problem in several function spaces. These results can be extended in a natural way to a broader class of elliptic systems with complex constant coefficients, as shown in \cite[Chapter 8]{MitMitMit-book5} using different techniques.

Given an open set $\Omega\subseteq{\mathbb{R}}^{n}$ and an aperture parameter
$\kappa\in(0,\infty)$, the {\tt nontangential approach region} at $x\in\partial\Omega$ is defined as
\begin{equation}\label{eq:approach-region}
\Gamma_{\kappa}(x):=\big\{y\in\Omega:\,|y-x|<(1+\kappa)\,{\rm dist}\,(y,\partial\Omega)\big\}, \qquad x\in\partial\Omega.
\end{equation}
From \cite[Proposition 2.9]{HofMitTay09}, we know that if $\Omega$ is an Ahlfors regular domain then $x\in\overline{\Gamma_\kappa(x)}$ for every $\kappa>0$ and $\mathcal{H}^{n-1}$-a.e. $x\in\partial\Omega$.
Next, for each continuous function $u$ on $\Omega$, we define the 
{\tt nontangential maximal operator} $\mathcal{N}_{\kappa}u:\partial\Omega\longrightarrow[0,+\infty]$ 
as 
\begin{equation}\label{eq:nt-max}
({\mathcal{N}}_{\kappa}u)(x):=
\sup_{y \in \Gamma_\kappa (x)} |u(y)|, \qquad x\in\partial\Omega,
\end{equation}
where we set $\big({\mathcal{N}}_{\kappa}u\big)(x):=0$ whenever
$\Gamma_\kappa(x)=\varnothing$. Additionally, we define the {\tt$\kappa$- nontangential limit} of $u$ at $x\in\partial\Omega$ as 
\begin{equation}
\big(u\big|^{{}^{\kappa-{\rm n.t.}}}_{\partial\Omega}\big)(x) :=
\lim_{\Gamma_\kappa(x) \ni y \to x} u(y),
\end{equation}
whenever this limit exists.

Set $\sigma:={\mathcal{H}}^{n-1}\lfloor\partial\Omega$ and
recall the harmonic double layer potential $K_\Delta$ defined in \eqref{eq:def-Kharmonic}. It is proved in \cite[Theorem 6.13]{HofMitTay09} that for every $p\in(1,\infty)$,
\begin{equation}
\frac{1}{2}I+K_{\Delta}:L^p(\partial\Omega,\sigma)\to L^p(\partial\Omega,\sigma)
\,\,\, \text{ is injective}
\end{equation}
whenever $\Omega$ is a bounded Ahlfors regular domain satisfying a two-sided local 
John condition such that $\mathbb{R}^n \setminus \overline{\Omega}$ is connected
and $\dist(\nu,[\VMO(\partial\Omega,\sigma)]^n)$ is small enough.
From this, Theorem \ref{thm:SIOCompact-unweighted}, and Proposition \ref{prop:Fred-S} we immediately obtain the following result.

\begin{proposition}\label{prop:invertibility}
Let $\Omega\subset\mathbb{R}^{n}$ be a bounded Ahlfors regular domain satisfying a two-sided local 
John condition such that $\mathbb{R}^n \setminus \overline{\Omega}$ is connected, and abbreviate 
$\sigma:={\mathcal{H}}^{n-1}\lfloor\partial\Omega$. 
Then, the following holds:

\begin{enumerate}
\item Fix $p\in(1,\infty)$ and a weight $w\in A_p(\partial\Omega)$. Then, there is $\delta>0$ such that if $\dist(\nu,[\VMO(\partial\Omega,\sigma)]^n)<\delta$,
\begin{equation}
\frac{1}{2}I+K_{\Delta}:L^p(\partial\Omega,w)\to L^p(\partial\Omega,w)
\,\,\, \text{ is invertible.}
\end{equation}
\item Let $\mathbb{X}$ be a RIBFS over $(\partial\Omega,\sigma)$ such that $1<p_{\mathbb{X}}\leq q_{\mathbb{X}} < \infty$ and fix a weight
$w \in A_{p_{\mathbb{X}}}(\partial \Omega)$ with $w(\partial\Omega)=\sigma(\partial\Omega)$.  Then, there exists $\delta>0$ such that if $\dist(\nu,[\VMO(\partial\Omega,\sigma)]^n)<\delta$,
\begin{equation}
\frac{1}{2}I+K_{\Delta}:\mathbb{X}(w)\to \mathbb{X}(w)
\,\,\, \text{ is invertible.}
\end{equation}
\item Let $p(\cdot)\in{\rm LH}(\partial\Omega)$ satisfy $1<p_{-} \leq p_+ < \infty$. Then, there is $\delta>0$ such that if $\dist(\nu,[\VMO(\partial\Omega,\sigma)]^n)<\delta$,
\begin{equation}
\frac{1}{2}I+K_{\Delta}:L^{p(\cdot)}(\partial\Omega,\sigma)\to L^{p(\cdot)}(\partial\Omega,\sigma)
\,\,\, \text{ is invertible.}
\end{equation}
\end{enumerate}
\end{proposition}

Next, fix $\Omega\subset\mathbb{R}^n$ as above and a function 
$f\in L^1(\partial\Omega\,,\, \sigma(x) / (1+|x|^{n-1}))$,
and define the action of the boundary-to-domain harmonic double layer potential operator 
${\mathcal{D}}_{\Delta}$ on $f$ as
\begin{equation}\label{eq:D}
{\mathcal{D}}_{\Delta}f(x):=\frac{1}{\omega_{n-1}}\int_{\partial\Omega}
\frac{\langle\nu(y),y-x\rangle}{|x-y|^{n}}f(y)\,d\sigma(y),\qquad\forall\,x\in\Omega,
\end{equation}
which verifies
\begin{equation}\label{eq:DProp}
{\mathcal{D}}_\Delta f\in{\mathscr{C}}^{\,\infty}(\Omega), \quad\,
\Delta({\mathcal{D}_\Delta} f)=0 \,\text{ in }\, \Omega, \quad\,
\mathcal{D}_{\Delta}f\big|_{\partial\Omega}^{{}^{\kappa-{\rm n.t.}}}
=\big(\tfrac{1}{2}I+K_{\Delta}\big)f \,\,\sigma\text{-a.e. on }\partial\Omega.
\end{equation}
Additionally, there is a non-decreasing function $\Phi: [1, \infty) \to [1, \infty)$ (depending on $n$, $p$, $\kappa$, and the geometric constants of $\partial\Omega$) such that
\begin{equation} \label{eq:Df-Lpw}
\|\mathcal{N}_{\kappa}(\mathcal{D}_{\Delta}f)\|_{L^p(\partial\Omega,w)} \leq \Phi([w]_{A_p}) \|f\|_{L^p(\partial\Omega,w)}
\qquad \forall f \in L^p(\partial\Omega,w)
\end{equation}
(cf., e.g., \cite[Chapter 3]{MarMarMitMitMit22}).
As a consequence, it is clear that whenever the operator $\frac{1}{2}I+K_\Delta$ is invertible, the function
\begin{equation}\label{eq:u-sol}
u(x) := \mathcal{D}_\Delta \Big( \Big(\frac{1}{2}I+K_\Delta \Big)^{-1} f\Big)(x),
\qquad x \in \Omega,
\end{equation}
satisfies $\Delta u=0$ in $\Omega$ and $u\big|_{\partial\Omega}^{{}^{\kappa-{\rm n.t.}}}=f\in L^p(\partial\Omega,w)$.
We close this section by combining these ideas to prove the following well-posedness result.

\begin{theorem}
Let $\Omega\subset\mathbb{R}^{n}$ be a bounded Ahlfors regular domain satisfying a two-sided local 
John condition such that $\mathbb{R}^n \setminus \overline{\Omega}$ is connected. Abbreviate 
$\sigma:={\mathcal{H}}^{n-1}\lfloor\partial\Omega$, and fix $\kappa\in(0,\infty)$. 
Then, the following holds:
\begin{enumerate}
\item Fix $p\in(1,\infty)$ and a weight $w\in A_p(\partial\Omega)$. Then, there is $\delta>0$ such that if $\dist(\nu,[\VMO(\partial\Omega,\sigma)]^n)<\delta$, the Dirichlet problem
\begin{equation}\label{eq:Dirichlet-Lpw}
\left\lbrace
\begin{array}{l}
u\in\mathscr{C}^{\infty}(\Omega),
\\[4pt]
\Delta u=0\,\,\text{ in }\,\,\Omega, 
\\[4pt] 
{\mathcal{N}}_{\kappa}u\in L^p(\partial\Omega,w),
\\[4pt]
u\big|_{\partial\Omega}^{{}^{\kappa-{\rm n.t.}}}=f\in L^p(\partial\Omega,w),
\\[4pt]
\end{array}
\right.
\end{equation}
is well-posed. Specifically, there is a unique solution $u$ given by \eqref{eq:u-sol}, and a constant $C\in(0,\infty)$ such that
\begin{equation}\label{eq:Dirichlet-esti-Lpw}
\|\mathcal{N}_\kappa u\|_{L^p(\partial\Omega,w)} \leq C \|f\|_{L^p(\partial\Omega,w)}.
\end{equation}
\item Let $\mathbb{X}$ be a RIBFS over $(\partial\Omega,\sigma)$ such that $1<p_{\mathbb{X}}\leq q_{\mathbb{X}} < \infty$ and fix a weight
$w \in A_{p_{\mathbb{X}}}(\partial \Omega)$ with $w(\partial\Omega)=\sigma(\partial\Omega)$.  Then, there exists $\delta>0$ such that if $\dist(\nu,[\VMO(\partial\Omega,\sigma)]^n)<\delta$, the Dirichlet problem
\begin{equation}\label{eq:Dirichlet-Xw}
\left\lbrace
\begin{array}{l}
u\in\mathscr{C}^{\infty}(\Omega),
\\[4pt]
\Delta u=0\,\,\text{ in }\,\,\Omega, 
\\[4pt] 
{\mathcal{N}}_{\kappa}u\in \mathbb{X}(w),
\\[4pt]
u\big|_{\partial\Omega}^{{}^{\kappa-{\rm n.t.}}}=f\in \mathbb{X}(w),
\\[4pt]
\end{array}
\right.
\end{equation}
is well-posed. Specifically, there is a unique solution $u$ given by \eqref{eq:u-sol}, and a constant $C\in(0,\infty)$ such that
\begin{equation}\label{eq:Dirichlet-esti-Xw}
\|\mathcal{N}_\kappa u\|_{\mathbb{X}(w)} \leq C \|f\|_{\mathbb{X}(w)}.
\end{equation}
\item Let $p(\cdot)\in{\rm LH}(\partial\Omega)$ satisfy $1<p_{-} \leq p_+ < \infty$. Then, there is $\delta>0$ such that if $\dist(\nu,[\VMO(\partial\Omega,\sigma)]^n)<\delta$, the Dirichlet problem
\begin{equation}\label{eq:Dirichlet-Lpvar}
\left\lbrace
\begin{array}{l}
u\in\mathscr{C}^{\infty}(\Omega),
\\[4pt]
\Delta u=0\,\,\text{ in }\,\,\Omega, 
\\[4pt] 
{\mathcal{N}}_{\kappa}u\in L^{p(\cdot)}(\partial\Omega,\sigma),
\\[4pt]
u\big|_{\partial\Omega}^{{}^{\kappa-{\rm n.t.}}}=f\in L^{p(\cdot)}(\partial\Omega,\sigma),
\\[4pt]
\end{array}
\right.
\end{equation}
is well-posed. Specifically, there is a unique solution $u$ given by \eqref{eq:u-sol}, and a constant $C\in(0,\infty)$ such that
\begin{equation}\label{eq:Dirichlet-esti-Lpvar}
\|\mathcal{N}_\kappa u\|_{L^{p(\cdot)}(\partial\Omega,\sigma)} \leq C \|f\|_{L^{p(\cdot)}(\partial\Omega,\sigma)}.
\end{equation}
\end{enumerate}
\end{theorem}

\begin{proof}
It follows from Proposition \ref{prop:invertibility} and \eqref{eq:DProp} that the function $u$ given by \eqref{eq:u-sol} is a solution for each of the Dirichlet problems. 
To establish uniqueness, recall that, under the hypotheses above, the spaces $L^p(\partial\Omega,w)$, $\mathbb{X}(w)$, and $L^{p(\cdot)}(\partial\Omega,\sigma)$ are all contained in $L^q(\partial\Omega,\sigma)$ for some $q\in(1,\infty)$ (as seen in the proof of Proposition \ref{prop:Fred-S}\ref{prop:Fred-S-d}). As a consequence, a solution for any of the three Dirichlet problems above also solves the corresponding Dirichlet problem in $L^q(\partial\Omega,\sigma)$, for which uniqueness is already known (cf. \cite[Theorem 7.2]{HofMitTay09}).
Finally, the estimates \eqref{eq:Dirichlet-esti-Lpw}, \eqref{eq:Dirichlet-esti-Xw}, and \eqref{eq:Dirichlet-esti-Lpvar} are consequences of \eqref{eq:Df-Lpw} and extrapolation (cf. Theorem \ref{thm:extrapolation-BFS}, along with Remarks \ref{rem:RdF}, \ref{rem:extrapolation-X}, and \ref{rem:extrapolation-Lpvar}).
\end{proof}


\def\ocirc#1{\ifmmode\setbox0=\hbox{$#1$}\dimen0=\ht0 \advance\dimen0
by1pt\rlap{\hbox to\wd0{\hss\raise\dimen0
\hbox{\hskip.2em$\scriptscriptstyle\circ$}\hss}}#1\else {\accent"17 #1}\fi}
\def\cprime{$'$} \def\cprime{$'$} \def\cprime{$'$}
\end{document}